\documentclass[final,a4paper]{elsarticle}
\usepackage[centering,top=1in,bottom=1in,left=1in,right=1in]{geometry}
\usepackage[format=hang]{caption}
\usepackage{subcaption}
\usepackage{amsmath}
\usepackage{physics}
\usepackage{graphicx}
\usepackage{mathtools}
\usepackage{comment}
\usepackage{placeins}
\usepackage{amsfonts}
\usepackage{tabularx}
\usepackage{lineno}
\usepackage{setspace}
\usepackage{xcolor}
\usepackage{bibspacing}
\usepackage{verbatim}

\usepackage{hyperref}

\hypersetup{
	colorlinks=true,
	linkcolor=black,
	citecolor=black,
	filecolor=black,
	urlcolor=black,
}

\makeatletter
\def\ps@pprintTitle{%
	\let\@oddhead\@empty
	\let\@evenhead\@empty
	\let\@oddfoot\@empty
	\let\@evenfoot\@oddfoot
}
\makeatother

\graphicspath{{./figures/}}

\bibstyle{elsarticle-num}
\biboptions{sort&compress}

\makeatletter
\def\ps@pprintTitle{%
	\let\@oddhead\@empty
	\let\@evenhead\@empty
	\let\@evenfoot\@oddfoot
}
\makeatother

\begin{document}
	
	\begin{frontmatter}
		
		\title{Optimal Structures for Failure Resistance Under Impact}
		
		\author{Andrew Akerson}
		
		\address{Division of Engineering and Applied Science, California Institute of Technology, Pasadena, CA 91125, USA\\}				

		\begin{abstract}
		The complex physics and numerous failure modes of structural impact creates challenges when designing for impact resistance. While simple geometries of layered material are conventional, advances in 3D printing and additive manufacturing techniques have now made tailored geometries or integrated multi-material structures achievable. Here, we apply gradient-based topology optimization to the design of such structures. We start by constructing a variational model of an elastic-plastic material enriched with gradient phase-field damage, and present a novel method to efficiently compute its transient dynamic time evolution. We consider a finite element discretization with explicit updates for the displacements. The damage field is solved through an augmented Lagrangian formulation, splitting the operator coupling between the nonlinearity and non-locality. Sensitivities over this trajectory are computed through the adjoint method, resulting in an adjoint problem which we solve in a similar manner to the forward dynamics. We demonstrate this formulation by studying the optimal design of 2D solid-void structures undergoing blast loading. Then, we explore the trade-offs between strength and toughness in the design of a spall-resistant structure composed of two materials of differing properties undergoing dynamic impact.
		\end{abstract}

	\end{frontmatter}

	\section{Introduction}

	The design of structures for impact or blast loading is encumbered by the complex interactions between wave propagation, plasticity, and material damage. This leads to failure modes such as plugging, fracture, petaling, and spall which are highly dependent on the material parameters, loading conditions, and structural layout~\cite{Backman1978}. This is further complicated by the trade-offs between properties such as strength and toughness when designing integrated structures of multiple materials. In practice, engineers typically start with industry standards and intuition, followed by sophisticated dynamical simulations to iterate on a design before it undergoes physical testing. Usually, these designs  consist of simple geometries of layered materials~\cite{Jena2009,Li2021,Huang2018}. However, with recent advances in additive manufacturing and 3D printing, we may now look to tailored designs with complex geometries and integrated materials~\cite{Rafiee2020,Ambrosi2020,Gadagi2020}. Additionally, the exponential growth of computational capabilities makes algorithmic optimal design methods feasible. This may allow us to efficiently design structures of unprecedented impact performance in scenarios where intuitive design is not sufficient. 
	
	Of the optimal structural design formulations, topology optimization has proven to be one of the most powerful methodologies. By considering the density of material at each point in the domain as a continuous variable, the design is posed as an optimization problem over these densities. Then, gradient-based optimization methods are used to iteratively update the design, where sensitivities are usually computed through the adjoint method. Originally introduced to optimize the compliance of linear elastic structures~\cite{Bendsoe1989}, topology optimization has since been applied to a wide range of applications including acoustic band-gaps~\cite{Sigmund2003}, piezoelectric transducers~\cite{Silva1999}, micro-electro-mechanical systems~\cite{Pedersen2002}, energy conversion devices~\cite{Collins2019}, and fluid structure interaction~\cite{Yoon2010}. 
	
	For the optimal design of impact problems, it is necessary to include transient dynamics, rate-dependent plasticity, and damage mechanics when modeling the material response. Past studies have addressed optimal design for transient dynamic evolution with elastic material models~\cite{Shobeiri2020,Nakshatrala2016}. Additionally, plasticity has been considered in both quasi-static~\cite{Schwarz2001,Wallin2016,Cox2019,Tauzowski2019} and dynamic settings~\cite{Nakshatrala2015,Ivarsson2018}. However, a structure with damage has only been considered in the static case. This has been studied in both the ductile~\cite{Li2017,Li2018} and quasi-brittle~\cite{Desai2022,Noel2017,Barbier2022} regime to design damage resistant structures. A variational mechanics model, where solutions are computed through energy principles, are favored to accurately model the physics and provide mathematical structure. Furthermore, an efficient computational method for these fields is necessary, as the iterative design process requires repeatedly simulating the dynamics for updated designs.
	
	To address the above mentioned requirements, we consider small-strain, rate-dependent plasticity enriched with continuum damage through a variational phase-field model in a transient dynamic setting. To efficiently simulate the dynamic response, we consider a finite element discretization where we employ an explicit update scheme for the displacement fields, and an implicit update for both the plasticity and damage. Because these irreversible damage updates are both nonlinear and non-local in nature, a direct computation would be prohibitively expensive. To this end, we use an operator-splitting augmented Lagrangian alternating direction method of multipliers. By introducing an auxiliary damage and Lagrange multiplier field, we accurately and efficiently solve the damage updates by iterating between a nonlinear local problem, a linear global problem, and a Lagrange multiplier update. 
	
	We look to optimize the material placement of the structure over the dynamic trajectory for a given objective function. By assuming the material parameters are dependent on a continuous design variable, we derive sensitivities through the adjoint method. This results in an adjoint dynamical system that we solve in a similar manner to the forward problem. These sensitivities are then used to update the design.
	
	We start in Section~\ref{sec:theoretical} by presenting the energy functional for system, then discuss the dynamic equilibrium relations. We apply the adjoint method, where sensitivities and adjoint relations are derived for a general objective. In Section~\ref{sec:numerics} we detail the solution process. First, we apply an augmented Lagrangian to operator split the damage updates. Then, using a finite element discretization, we solve the system with explicit displacement updates, followed by implicit plasticity and damage updates. We demonstrate the accuracy and efficiency of the numerical scheme by considering the solution convergence and time-scaling for a model problem. We use a similar numerical scheme for the adjoint system and the associated dual variables. Next, in Section~\ref{sec:interp}, we discuss material interpolation schemes through intermediate densities for both solid-void structures and multi-material designs.  In Section~\ref{sec:results}, we demonstrate the methodology by looking at two examples. First we consider the design of 2D solid-void structures optimized for blast loading. Next, we explore the trade-offs between strength and toughness in a two material spall-resistant structure undergoing impact. Finally, in Section~\ref{sec:conclusions}, we summarize our findings and discuss further directions.

	\section{Theoretical Formulation} \label{sec:theoretical}
	
	\subsection{Forward Problem} \label{sec:forward}

	We consider an elastic-plastic material capable of sustaining damage occupying a bounded, open domain $\Omega \subset \mathbb{R}^n$ in its reference configuration over time $[0, T]$.  We assume prescribed loads on $\partial_f \Omega \subset \partial \Omega$ and prescribed displacements on   $\partial_u \Omega \subset \partial \Omega$. We consider small-strain, rate-dependent J-2 plasticity with isotropic hardening to model the plasticity~\cite{Ortiz1989,lubliner2008plasticity}. Damage is measured by the phase-field scalar quantity $a : \Omega \times [0, T] \mapsto [0, 1]$, where values of $0$ and $1$ correspond to the undamaged and fully damaged states. Here, we use a phase-field fracture model which we adapt for damage by considering a finite length scale~\cite{Bourdin2000}. These models have been modified for ductile fracture by including small-strain plasticity~\cite{Brach2019}, and we adopt a similar formulation. We assume the material parameters are dependent on a design field $\eta: \Omega \times [0, T] \mapsto [0, 1]$ which determines the species of material at each point. We consider a variational structure, where minimization principles yields the internal variable evolution~\cite{Ortiz1999}. Thus, we consider the incremental energy
	\begin{equation}
	\begin{aligned}
		\mathcal{E} (u, q, \varepsilon^p, a, \eta) = \int_{\Omega} \Bigg{\{}  &W^e(\varepsilon, \varepsilon^p, a, \eta) +   d(a) \left[W^p(q, \eta) + \int_0^t g^*(\dot{q}, \eta) \, dt \right] \\
		&+ \frac{G_c(\eta)}{4 c_w} \left[ \frac{w^a(a, \eta)}{\ell(\eta)} + \ell(\eta) \norm{\nabla a}^2 \right]  + \int_{0}^t \psi^*(\dot{a}, \eta) dt  \Bigg{\}} d\Omega,
	\end{aligned}
	\end{equation}
	where $u : \Omega \times [0, T] \mapsto \mathbb{R}^n$ is the displacement field, $\varepsilon^p : \Omega \times [0, T] \mapsto \mathbb{R}^{n\times n}$ is the volume preserving plastic strain, and $q : \Omega \times [0, T] \mapsto \mathbb{R}^+$ is the accumulated plastic strain whose evolution is defined by 
	\begin{equation}
    \dot{q} = \sqrt{\frac{2}{3} \dot{\varepsilon}^p \cdot \dot{\varepsilon}^p}. 
    \end{equation}
    $W^e$ is the stored elastic energy density, which accounts for the tension-compression asymmetry in its damage dependence~\cite{Amor2009},
    \begin{equation}
        W^e(\varepsilon, \varepsilon^p, a, \eta) = \frac{K(\eta)}{2} \tr^-(\varepsilon^e)^2 + d(a)\left[ \frac{K(\eta)}{2} \tr^+(\varepsilon^e)^2 + \mu (\eta) \varepsilon^e_D : \varepsilon^e_D \right],
    \end{equation}
    where $K$ and $\mu$ are the bulk and shear moduli. $d(a)$ models the weakening of the material with damage,
    \begin{equation}
        d(a) = (1 - a)^2 + d_1 a^2,
    \end{equation}
   where $d_1 << 1$. $\varepsilon^e = \varepsilon - \varepsilon^p$ is the elastic strain, and $\varepsilon^e_D$ is its deviatoric component.  $\tr^+(\varepsilon)$ and $\tr^-(\varepsilon)$ are the positive and negative parts of the strain trace, 
    \begin{equation}
        \tr^+(\varepsilon) = \max(\tr(\varepsilon), 0), \qquad \tr^-(\varepsilon) = \min(\tr(\varepsilon), 0).
    \end{equation}
    This decomposition of the volumetric strain allows for tension-compression asymmetry in the damage model; the tensile bulk modulus is affected by damage, while the compressive bulk modulus remains unaffected. $W^p$ and $w^a$ are the plastic and damage hardening functions, respectively. The damage parameters $G_c$ and $\ell$ control the toughness and damage length scale, with $c_w$ as a normalization constant. Finally, the rate dependence of both the damage and plastic hardening is handled by the dissipation potentials $\psi^*$ and $g^*$, respectively. These functions also account for irreversibility, as they take a value of $+\infty$ for negative rates,
    \begin{equation}
        g^*(\dot{q}) = \begin{cases}
        \bar{g}^*(\dot{q}) \quad & \dot{q} \geq 0 \\
        \infty \quad & \dot{q} < 0
        \end{cases}, \qquad 
         \psi^*(\dot{a}) = \begin{cases}
        \bar{\psi}^*(\dot{a}) \quad & \dot{a} \geq 0 \\
        \infty \quad & \dot{a} < 0
        \end{cases}.
    \end{equation}
    For the plastic potentials, we consider power-law hardening and rate-sensitivity functions
    \begin{equation}
        W^p(q) = \sigma_y\left[ q + \frac{n \varepsilon^p_0}{n + 1}\left( \frac{q}{\varepsilon^p_0}\right)^{(n+1)/n}\right], \quad \bar{g}^*(\dot{q}) = \frac{ m \sigma_y \dot{\varepsilon}^p_0}{m + 1}\left( \frac{\dot{q}}{\dot{\varepsilon}^p_0}\right)^{(m + 1)/m}.
    \end{equation}
    $\varepsilon^p_0$ and $\dot{\varepsilon}^p_0$ are the reference plastic strain and strain rate and $\sigma_y$ is the initial yield stress.  $n$ and $m$ are the powers for the hardening and rate sensitivity, with the perfecty plastic and rate-indepdendent cases occurring as $n \rightarrow \infty+$ and $m \rightarrow \infty+$, respectively~\cite{Ortiz1989}. For the damage hardening, we consider a quadratic function
    \begin{equation}
        w^a(a) = w_1 a + (1 - w_1) a^2,
    \end{equation}
    where $w_1 \in [0, 1]$, which ensures $w^a(1) = 1$. For simplicity, we consider the damage to be rate-independent by choosing $\bar{\psi}^*(\dot{a}) = 0$.  Here, we scale both the plastic potential and shear modulus with the same damage function $d(a)$. Thus, the yield strength and Mises stress have the same damage dependence, leading to damage independent plastic updates. 
    
    We consider dynamic evolution through the incremental action integral
	\begin{equation}
		\begin{aligned}
			\mathcal{L}(u, q, \varepsilon^p, a, \eta) = & \int_{t_1}^{t_2} \left \{ \mathcal{E}(u, q, \varepsilon^p, a, \eta) - \int_{\Omega}  \frac{\rho(\eta)}{2} \abs{\dot{u}}^2 d \Omega  - \int_{\Omega} f_b \cdot u \, d\Omega -  \int_{\partial_f \Omega} f \cdot u \, dS  \right \} \, dt, \\
		\end{aligned}
	\end{equation}
	where $f_b$ and $f$ are the body force and surface tractions, and $\rho$ is the material density. Stationarity of this action integral gives the dynamic evolution and the kinetics of the internal variables~\cite{Mielke2005}
	\begin{subequations} \label{eq:forward_eqs_continuous}
	\begin{align} 
			0 &= \int_{\Omega} \left[ \rho \ddot{u} \cdot \delta u +  \pdv{W^e}{\varepsilon}  \cdot \nabla \delta  u  \right] \, d\Omega - \int_\Omega f_b \cdot \delta u \, d\Omega - \int_{\partial_2\Omega} f \cdot \delta u \, d\Omega  \quad &&\forall \delta u \in \mathcal{U}, \label{eq:for_u} \\
			\qquad 0 &\in \bar{\sigma}_M - \pdv{W^p}{q} - \partial g^*,  &&  \text{ on } \Omega, \label{eq:for_yield}\\
			\qquad 0 &= \dot{\varepsilon^p} - \dot{q} M  &&  \text{ on } \Omega, \label{eq:for_flow}\\
			0 &\in \pdv{W^e}{a} + \pdv{d}{a} \left(W^p + \int_0^t g^*(\dot{q}) \, dt \right) - \nabla \cdot \left ( \frac{G_c \ell }{2 c_w} \nabla a \right) + \frac{G_c}{4 c_w \ell} \pdv{w^a}{a} + \partial \psi^*  && \text{ on } \Omega, \label{eq:for_a}\\
			& a = 0 \ \text{ on } \partial_u \Omega, \qquad \nabla a \cdot n = 0 \ \text{ on } \partial_f \Omega. \label{eq:for_a_bc}
	\end{align}
	\end{subequations}
	Here, we assume quiescent initial conditions. $\mathcal{U}$ is the space of admissible displacement variations
	\begin{equation}
	    \mathcal{U} = \{ u \in H^1(\Omega), \ u = 0 \text{ on } \partial_u \Omega\}.
	\end{equation}\eqref{eq:for_u} is the second-order dynamic evolution of the displacement field. \eqref{eq:for_yield} and \eqref{eq:for_flow} are the yield relation and the evolution of the plastic strain, where $\bar{\sigma}_M$ is the normalized Mises stress (divided through by d(a)), and $M$ is the direction of plastic flow. \eqref{eq:for_a} is the irreversible evolution of the damage field, with \eqref{eq:for_a_bc} being the boundary conditions for $a$. The differential inclusion in the yield relation and damage equilibrium enforces the irreversibility of their respective internal variables.

	\subsection{ Sensitivities and Adjoint Problem} \label{sec:adjoint} 
	We look to find the design field $\eta(x)$ such that an objective, dependent on the dynamic trajectory, is minimized. Thus, we consider a general objective of integral form
	\begin{equation}
		\begin{aligned}
			& \underset{\eta(x)}{\min} \quad \mathcal{O}(\eta) := \int_0^T  \int_{\Omega} o(u, q, \varepsilon^p, a, \eta) \, d\Omega \ dt \\
			& \text{subject to:}  \qquad \text{Equillibrium relations in \eqref{eq:forward_eqs_continuous}.}
		\end{aligned}
	\end{equation} 
	To conduct gradient-based optimization, the variation of the objective with $\eta$ must be computed. For this, we employ the adjoint method~\cite{Plessix2006}. We introduce fields $\xi$, $\gamma$, $\mu$, and $b$ as the dual variables to the displacement, plastic hardening, plastic strain, and the damage fields, respectively. We consider the necessary Kuhn-Tucker conditions for the irreversible equilibrium relations, and carry out the adjoint calculation. The full details of this can be found in \ref{app:adjsens}. This gives the total variation of the objective as
	\begin{equation} \label{eq:sensitivities}
		\begin{aligned}
			\mathcal{O}_{, \eta} \delta \eta = &\int_{0}^T  \int_{\Omega} \bigg\{ \pdv{o}{\eta} + \pdv{\rho}{\eta} \ddot{u} \cdot \xi + \pdv{^2 W^e}{\varepsilon \partial \eta} \cdot \nabla \xi + b \dot{a} \left( \pdv{^2 W^e}{a \partial \eta} + \pdv{d}{a} \pdv{W^p}{\eta} + \pdv{d}{a}  \int_{0}^t \pdv{g^*}{\eta} d \tau \right) \\
			& + \frac{1}{2 c_w}\pdv{(G_c \ell)}{\eta} \nabla(b \dot{a}) \cdot \nabla{a} + b \dot{a} \left( \frac{w^{a \prime}}{4 c_w}\pdv{(G_c/ \ell)}{\eta} + \pdv{^2 {\psi}^*}{\dot{a} \partial \eta} \right) + \gamma \dot{q} \left( \pdv{\bar{\sigma}_M}{\eta} - \pdv{\sigma_0}{\eta} - \pdv{^2 {g}^*}{\dot{q} \partial \eta}\right) \bigg\} \ \delta \eta \, d\Omega \, dt,
		\end{aligned}
	\end{equation}
	where the adjoint variables satisfy the dynamic evolution 
	\begin{subequations} \label{eq:adjoint_system_continuous}
    \begin{align}
    &0 = \int_{\Omega} \left[ \rho \ddot{\xi} \cdot \delta_\eta u + \pdv{o}{u} \cdot \delta_\eta u  + \left( \nabla \xi \cdot \pdv{^2 W^e}{\varepsilon \partial \varepsilon}  + b \dot{a} \pdv{^2 W^e}{a \partial \varepsilon} +  \gamma \dot{q} \pdv{\bar{\sigma}_M}{\varepsilon} - \dot{q} \mu \cdot \pdv{M}{\varepsilon} \right) \cdot \nabla \delta_\eta u \right] \ d\Omega && \forall \delta_\eta u \in \mathcal{U} \label{eq:adju}\\
    & \dv{}{t} \left[ \gamma \left( \bar{\sigma}_M - \sigma_0 - \pdv{\bar{\psi}^*}{\dot{q}}\right)- \gamma \dot{q} \pdv{^2 \bar{g}^*}{\dot{q}^2} + \pdv{\bar{g}^*}{\dot{q}} \left( \int_t^T b \dot{a} d^{\prime}(a) d\tau\right) - \mu \cdot M \right]  && \nonumber \\
    & \hspace{5cm} =\pdv{o}{q} + b\dot{a} d^{\prime}(a) \pdv{W^p}{q} - \gamma \dot{q}  \pdv{\sigma_0}{q} && \text{ on } \Omega \label{eq:adjp}\\
    & \dv{\mu}{t} = \pdv{o}{\varepsilon^p} + \nabla \xi \cdot \pdv{^2 W^e}{\varepsilon \partial \varepsilon^p} + b \dot{a} \pdv{^2 W^e}{a \partial \varepsilon^p} + \gamma \dot{q} \pdv{\bar{\sigma}_M}{\varepsilon^p} - \dot{q} \mu \cdot \pdv{M}{\varepsilon^p} && \text{ on } \Omega  \label{eq:adjpstrain}\\
    &\dv{}{t} \left[ D_a b + \pdv{^2 \bar{\psi}^* }{\dot{a}^2} b \dot{a}\right] =  \pdv{o}{a} + \pdv{^2 W^e}{a \partial \varepsilon} \cdot \nabla \xi + b \dot{a} \left( \pdv{^2 W^e}{a^2} + \frac{G_c}{4 c_w \ell} \pdv{^2 w^a}{a^2}\right) && \nonumber \\
    & \hspace{5cm} + b \dot{a} d^{\prime \prime} \left(  W^p + \int_0^t g^* d\tau \right) - \nabla \cdot \left( \frac{G_c \ell}{2 c_w} \nabla(b \dot{a}) \right) && \text{ on } \Omega \label{eq:adja}\\
    & \quad \xi |_{t = T} = 0, \quad  \dot{\xi} |_{t = T} = 0, \quad \gamma |_{t = T} = 0, \quad  \mu |_{t = T} = 0, \quad b |_{t = T} = 0, \nonumber
    \end{align}
    \end{subequations}
	 where
	\begin{equation}
    D_a = \pdv{W^e}{a}  + \pdv{d}{a} \left(W^p + \int_0^t g^*\, d\tau \right) - \nabla \cdot \left( \frac{G_c \ell }{2 c_w} \nabla a \right) + \frac{G_c}{4 \ell c_w }\pdv{w^a}{a}  + \pdv{\bar{\psi}^*}{\dot{a}}.
    \end{equation}
    These are dependent on the forward problem solution and must be solved backwards in time. Once the forward problem is solved in time for $u(t), a(t), q(t),$ and  $\varepsilon^p(t)$, they can be used to solve the adjoint problem backwards in time for $\xi(t), b(t), \gamma(t),$ and $\mu(t)$. The sensitivities can then be computed from \eqref{eq:sensitivities}. Details of the numerical methods to solve the forward and adjoint problem are discussed in the proceeding section.
	
	\section{Numerics} \label{sec:numerics}
	\subsection{Forward Problem}
	We discuss the details for the numerical evolution of the forward dynamics. First we introduce an augmented Lagrangian formulation to split the nonlinear and non-local operator coupling in the damage field equilibrium. Then, using a finite element discretization, we discuss the computational procedure for updating the displacements, plasticity and damage variables. Finally, we study the accuracy and efficiency of our formulation by studying the solution behavior for varying mesh sizes. 
	\subsubsection{ Augmented Lagrangian} \label{sec:auglag}
	The differential inclusion and gradient terms in the damage evolution of \eqref{eq:for_a} result in a nonlinear and non-local state equation for the damage updates. While there exist methods to directly solve these non-local constrained problems, they result in expensive computations that would be required at every timestep. Thus, we consider an augmented Lagrangian formulation to split this operator, and solve the system using an alternating direction method of multipliers (ADMM)~\cite{glowinski1989augmented,Fortin1983AugmentedLM}. This method has been used to efficiently solve non-linear elasticity problems with internal variable evolution~\cite{Zhou2021}. We introduce the auxiliary field $\alpha \in L^2(\Omega)$ and constrain $a = \alpha$ weakly for all time with the Lagrange multiplier $\lambda \in L^2(\Omega)$ and penalty factor $r$. Thus, we consider the modified incremental energy
	\begin{equation}
	\begin{aligned}
		\mathcal{E}  = \int_{\Omega} \Bigg{\{}  &W^e(\varepsilon, \varepsilon^p, \alpha, \eta) +   d(\alpha) \left[W^p(q, \eta) + \int_0^t g^*(\dot{q}, \eta) \, dt \right] \\
		&+ \frac{G_c(\eta)}{4 c_w} \left[ \frac{w^a(\alpha, \eta)}{\ell(\eta)} + \ell(\eta) \norm{\nabla a}^2 \right]  + \int_{0}^t \psi^*(\dot{\alpha}, \eta) dt  + \frac{r}{2}(a - \alpha)^2 + \lambda(a - \alpha) \Bigg{\}} \ d\Omega.
	\end{aligned}
	\end{equation}
	Stationarity of the action integral using this augmented energy results in the equilibrium relations identical to that of \eqref{eq:forward_eqs_continuous}, with the exception that \eqref{eq:for_a} be replaced by
    \begin{subequations} \label{eq:admm}
	\begin{align} 
		&\lambda + r(a - \alpha)  - \pdv{W^e}{\alpha}  - d^{\prime}(\alpha) \left[W^p(q) + \int_0^t g^*(\dot{q}) \, dt \right] - \frac{G_c}{4 c_w \ell }\pdv{w^a}{\alpha} (\alpha)  \in \partial \psi^* \left( \dot{\alpha} \right) \quad && \text{ on } \Omega, \label{eq:admm1}\\
		&0 = \int_\Omega \left[ \frac{G_c \ell}{2 c_w} \nabla a \cdot \nabla \delta a + r (a - \alpha) \delta a + \lambda \delta a \right] \, d \Omega \quad && \forall \delta a \in \mathcal{A},  \label{eq:admm2} \\
		&0 = \int_{\Omega} (a - \alpha) \delta \lambda \, d \Omega  \quad && \forall \delta \lambda \in L^2(\Omega).  \label{eq:admm3} 
	\end{align}
	\end{subequations}
	where
	\begin{equation}
	    \mathcal{A} = \{a \in H^1(\Omega), \ a = 0 \text{ on } \partial_u \Omega \}.
	\end{equation}
	With $\alpha$ as the unknown, \eqref{eq:admm1} is a nonlinear local problem. Correspondingly, the second line \eqref{eq:admm2} is a linear global problem for $a$. The de-coupling of nonlinearity and non-locality allows for the efficient computation of the damage evolution, which we discuss with the numerical implementation.

	\subsubsection{Discretization and Solution Procedure} \label{sec:fem}
	We discretize the system with standard $p = 1$ Lagrange finite elements for the displacement field $u$ and the damage field $a$ as
    \begin{equation}
    u = \sum_{i = 1}^{n_u} u_i N^u_i(x), \qquad a = \sum_{i = 1}^{n_a} a_i N^a_i(x),
    \end{equation} 
    where $N^u_i \in \mathbb{R}^n$ and $N^a_i \in \mathbb{R}$ are standard vector and scalar valued first-order shape functions with compact support. The fields $\alpha$, $q$, and $\varepsilon^p$ are discretized at quadrature points
    \begin{equation}
    \alpha(x_g) = \alpha_g, \qquad q(x_g) = q_g, \qquad \varepsilon^p (x_g) = \varepsilon^p_g, 
    \end{equation}
    for some Gauss point $x_g$. The Lagrange multiplier field $\lambda$ is discretized in the same finite element space we use for $a$ as
    \begin{equation}
    \lambda =  \sum_{i = 1}^{n_a} \lambda_i N^a_i(x).
    \end{equation}
	Finally, the design field $\eta$ is assumed constant on each element. 
	
	We start with an explicit central difference scheme to update the displacement field. Because the plasticity updates do not depend on the damage field, $q$ and $\varepsilon^p$ are next computed implicitly with a backwards Euler update.  Finally, the damage field is updated implicitly by iterating between the nonlinear local problem for $\alpha$ by solving \eqref{eq:admm1}, the linear global problem for $a$ through \eqref{eq:admm2}, and a Lagrange multiplier update for $\lambda$ until convergence. Since the operator for the global problem remains identical between iterations, we need only construct the system matrix and perform the sparse LU decomposition once, where subsequent solves involve only a right-hand side assembly and back-substitution. For the $n$ to $n+1$ time-step the displacement updates are
    \begin{equation}
    \begin{aligned}
    \ddot{u}_i^{n} &= M^{-1}_{ij} F_j^n(u^n, \varepsilon^{p,n}, \alpha^n, t^n),\\
    \dot{u}_i^{n + 1/2} &= \dot{u}_i^{n - 1/2} + \Delta t^n \, \ddot{\bar{u}}_i^{n}, \\
    u_i^{n + 1} &= u_i^n + \Delta t^{n + 1/2} \, \dot{u}_i^{n + 1/2},
    \end{aligned}
    \end{equation}
    where,
    \begin{equation}
    M_{ij} = \int_{\Omega} \rho(x) N_i^u \cdot N_i^u \, d\Omega, \qquad F_j^n = \int_{\Omega} \left[ -\pdv{W^e}{\varepsilon} (\varepsilon^n, \varepsilon^{p,n},\alpha^n, \eta) \cdot \nabla N^u_j  + f_b \cdot N^u_j \right] \, d\Omega - \int_{\partial_f \Omega} f \cdot N^u_j \, d\Omega.
    \end{equation}
    In standard fashion, these integrals are approximated with Gauss quadrature. Again, since the plastic evolution does not depend on the damage field, we update the plasticity variables through an implicit backwards Euler discretization. For this, we employ a predictor-corrector scheme~\cite{Ortiz1989} to solve point-wise at each quadrature point,
    \begin{equation}
    \begin{aligned}
    & 0 \in \bar{\sigma}_M (\varepsilon^{n + 1} |_{x_g}, \varepsilon^{p, (n + 1)}_g, \eta(x_g) ) - \sigma_0(q^{n+1}_g, \eta(x_g) ) - \partial g^* \left( \frac{q^{n + 1}_g - q^n_g}{\Delta t}, \eta(x_g) \right), \\
    & \varepsilon^{p, (n+1)}_g = \varepsilon^{p, n}_g + \Delta q M (\varepsilon^{n + 1}_g, \varepsilon^{p, (n + 1)}_g).
    \end{aligned}
    \end{equation} 
    The update for $\alpha$  uses an implicit backwards Euler method, coupled with ADMM for the fields $a$ and $\lambda$. This reduces to iterations between a nonlinear point-wise problem for the updates of $\alpha$, a linear global problem for $a$, and an update for $\lambda$. 
    
    We summarize these operations for the $n$ to $n+1$ time-step. Given $u^{n+1}$, $q^{n+1}$, $\varepsilon^{p, (n+1)}$, we initialize values $\tilde{\lambda}^0 = \lambda^n$, $\tilde{a}^0 = a^n$, and iterate over $i$:
    \begin{itemize}
    	\item {\textit{Step 1: Non-linear local problem.}} Update $\tilde{\alpha}^{i+1}$ by solving at each $x_g$
    	\begin{equation}
    	\begin{aligned}
    	-\pdv{W^e}{\alpha} &\left(\varepsilon^{n+1}|_{x_g}, \tilde{\alpha}^{i + 1}_g, \eta(x_g) \right) -  d^{\prime}(\tilde{\alpha}^{i+1}_g) \left[W^p(q^{n+1}_g , \eta(x_g) ) + \int_0^t g^*(\dot{q}_g, \eta(x_g) ) \, dt \right] \\
    	&- \frac{G_c(\eta(x_g))}{4 c_w \ell(\eta(x_g))}\pdv{w^a}{\alpha} \left( \tilde{\alpha}^{i + 1}_g, \eta(x_g) \right) + \tilde{\lambda}^i|_{x_q} + r \left(\tilde{a}^i |_{x_g} - \tilde{\alpha}^{i + 1}_g \right) \in \partial \psi^* \left(\frac{\tilde{\alpha}^{i+1}_q - \alpha^{n}_q}{\Delta t_n}, \eta(x_q) \right).
    	\end{aligned}
    	\end{equation}
    	
    	\item {\textit{Step 2: Linear global problem.}} Update $\tilde{a}^{i+1}$ by solving
    	\begin{equation}
    	{K}_{pj} \, \tilde{a}^{i+1}_j = V_p(\tilde{\alpha}^{i+1}, \tilde{\lambda}^i),
    	\end{equation}
    	where
    	\begin{equation}
    	K_{pq} = \int_{\Omega} \left[ \frac{G_c (\eta) \ell (\eta)}{2 c_w} \nabla N^a_p \cdot  \nabla N^a_q + r N^a_p N^a_q \right] \, d\Omega, \qquad V_p(\alpha, \lambda) = \int_{\Omega} \left( r \alpha - \lambda \right) N^a_p \, d\Omega.
    	\end{equation}
    	
    	\item {\textit{Step 3: Update Lagrange multiplier.}} Update $\tilde{\lambda}^{i+1}$ by
    	\begin{equation}
    	\tilde{\lambda}^{i + 1}_j = \tilde{\lambda}^i_j + r(\tilde{a}^{i + 1}_j - S^{-1}_{jk} \hat{\alpha}^{i + 1}_k ),
    	\end{equation}
    	where
    	\begin{equation}
    	S_{jk} = \int_{\Omega} N^a_j N^a_k \, d\Omega, \qquad \hat{\alpha}^{i+1}_k = \int_{\Omega} \tilde{\alpha}^{i+1} N^a_k \, d\Omega.
    	\end{equation}
    	Note: this is the weak form of the update $\Delta \lambda = r(a - \alpha)$.
    	
    	\item {\textit{Step 4: Check for convergence.}} Check both primal and dual feasibility
    	\begin{equation}
    	\begin{aligned}
    	r_p &:= \norm{\bar{a}^{i+1} - \hat{\alpha}^{i+1}}_{l^2} \leq \frac{1}{\sqrt{n_a}} r^{tol}_{abs} + r^{tol}_{rel} \max \left(\norm{\hat{\alpha}^{i+1}}_{l^2}, \norm{\bar{a}^{i+1}}_{l^2} \right),    \\
    	r_d &:= r \norm{\bar{a}^{i+1} - \bar{a}^{i}}_{l^2} \leq \frac{1}{\sqrt{n_a}} r^{tol}_{abs} + r^{tol}_{rel} \norm{\bar{\lambda}^{i + 1}} ,
    	\end{aligned}
    	\end{equation} 
    	where
    	\begin{equation}
    	\bar{a}^{i + 1}_{j} = S_{jk} \tilde{a}^{i + 1}_k, \qquad \bar{\lambda}^{i+1}_j = S_{jk} \tilde{\lambda}_j.
    	\end{equation}
    	In the above, we use the vector $l^2$ norm
    	\begin{equation}
    	\norm{\bar{a}}^2_{l^2} = \sum_{i = 1}^{n_a} \bar{a}_i^2
    	\end{equation}
    \end{itemize}
    until convergence, and update $\alpha^{n+1} = \tilde{\alpha}^i$, $a^{n+1} = \tilde{a}^i$, and $\lambda^{n+1} = \tilde{\lambda}^i$. 
    For faster convergence, we update the penalty value $r$ between iterations. As larger values of $r$ improve primal feasibility convergence while slowing the dual feasibility convergence (and vice-versa), adapting the value of $r$ based on these feasibility values can lead to few iterations~\cite{Zhou2021,Boyd2010}. Thus, we consider the following scheme
    \begin{equation} \label{eq:updater}
        r = \begin{cases}
            \min(\gamma_r r, r_{max}) & \text{ if } r_p / r_d > \tau \\
             \max(r/\gamma_r , r_{min}) & \text{ if } r_d / r_p > \tau \\
             r & \text{ else }
        \end{cases}.
    \end{equation}
    In our study, we choose $\tau = 10$, and take $\gamma_r = 2$.
    
    \subsubsection{Accuracy and Efficiency}
    
    	\begin{figure}[]
		\begin{center}
				\includegraphics[width=0.95\textwidth]{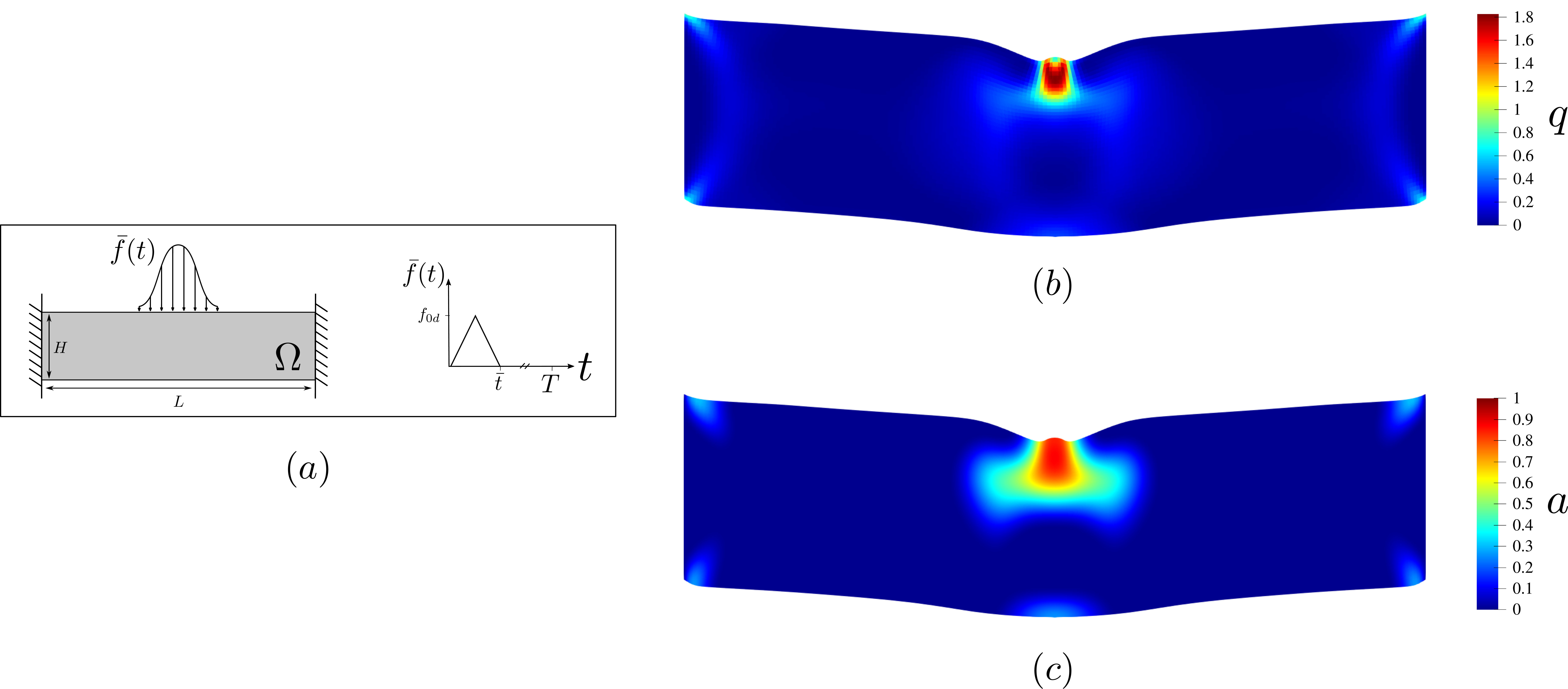}
		\end{center}
		\caption{The model problem we use to study the accuracy and efficiency of our formulation. We consider a rectangular geometry with a impulse Gaussian loading profile (a). Additionally, deformed configurations with accumulated plasticity (b) and damage fields (c) are shown at the final time-step computed on a 200$\times$50 mesh. }
		\label{fig:convdef}
	\end{figure}
    
    To analyze the efficiency and efficacy of the above formulation, we study a model problem. We consider a clamped bar undergoing dynamic loading on its top surface, as shown in Figure~\ref{fig:convdef}\hyperref[fig:convdef]{a}. The loading is chosen such that the structure undergoes both plastic and damage evolution along its trajectory. Figure~\ref{fig:convdef}\hyperref[fig:convdef]{b} and \ref{fig:convdef}\hyperref[fig:convdef]{c} show the plasticity and damage fields at the final time.  We investigate the solution convergence and time-scaling for uniform meshes varying from 60$\times$15 to 600$\times$150 for a constant 18,000 time-steps. Each of the simulations are run on 6 CPU cores using shared memory. The absolute and relative ADMM tolerance is set to a constant $r_{abs}^{tol} = r_{rel}^{tol} = 10^{-7}$.
    
\begin{figure}
    \centering
    \begin{subfigure}{0.5\textwidth}
      \centering
      \includegraphics[width=0.9\textwidth]{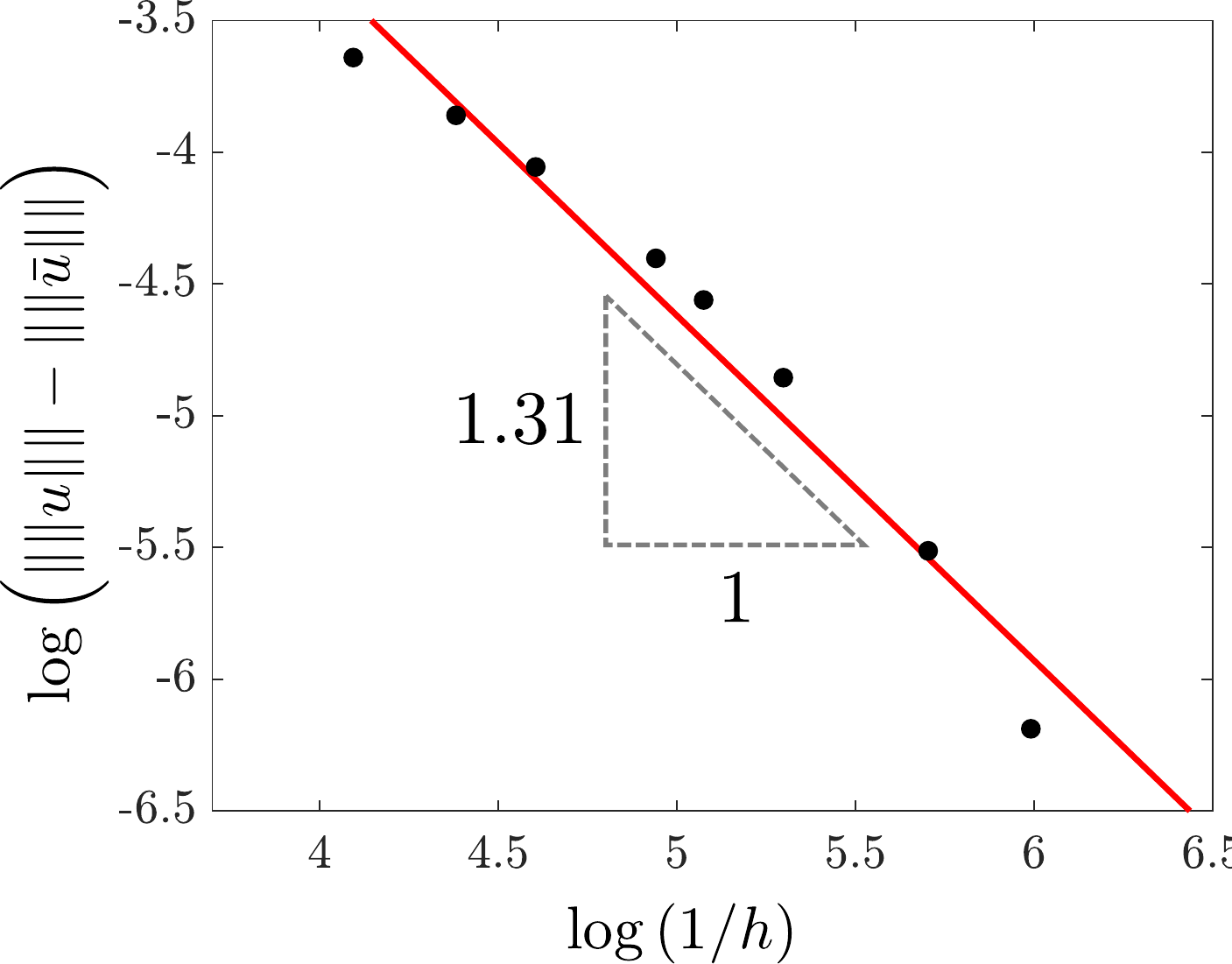}
      \subcaption{}
      \label{subfig:sub1conv}
    \end{subfigure}%
    \begin{subfigure}{0.5\textwidth}
      \centering
      \includegraphics[width=0.9\textwidth]{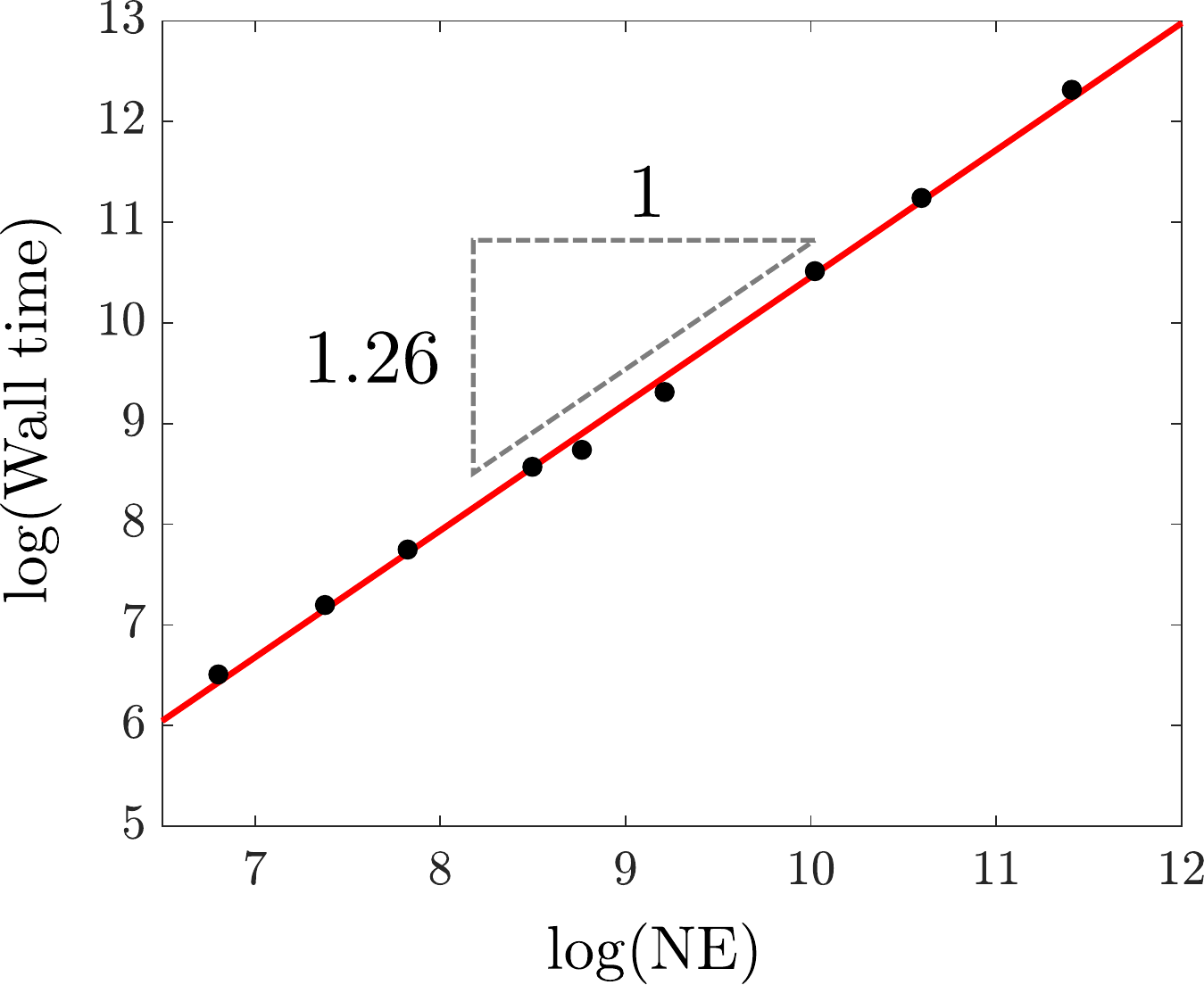}
      \subcaption{}
      \label{subfig:sub2conv}
    \end{subfigure}
    \caption{Solution convergence and time-scaling plots for varying mesh sizes. The solution norm $\| \| u \| \|$ is studied relative to the characteristic mesh size $h$ (a). For time-scaling, we consider the wall time v.s. the number of element, NE (b). The black dots represent data for each of the simulations, while the red lines show the linear fits, with the first order coefficients denoted on the triangles. }
    \label{fig:convergence}
    \end{figure}

    To study the solution convergence, we consider the $L^2$ norm in time of the $H^1$ norm in space, which we denote as $\| \| \cdot \| \| := \| \left( \| \cdot \|_{H^1(\Omega)} \right) \|_{L^2(0, T)} $. We investigate $\| \| u \| \| $ for the varying meshes. As an analytical solution does not exist, we consider the solution on the $600\times150$ mesh as the reference, $\bar{u}$. Figure~\ref{subfig:sub1conv} shows the convergence of the displacement norm for varying characteristic mesh size $h$. A linear fit yields a convergence rate of $1.31$, demonstrating super-linear convergence even while undergoing large plastic and damage evolution. Next, we study the time-scaling for varying mesh sizes. For meshes varying from 900 to 90,000 elements, we see a growth rate with wall time of $1.26$. This exceptional scaling may be attributed to the ADMM algorithm for computing the damage evolution. As the linear global problem has a constant operator for each penalty value $r$, these matrices may be pre-computed and treated with an LU decomposition in set-up. Then, each of the linear solves may be executed through efficient back-substitution. It is expected that this scaling  breaks down if the number of elements increases significantly, as the solution time is then dominated by the more inefficient LU decomposition. 
    
	\subsection{Adjoint Problem}
	We now turn to the details of the numerical evolution of the adjoint problem, which must be solved backwards in time using the solution to the forward problem. For efficiency, we employ another augmented Lagrangian formulation for the adjoint damage variable update. Then, we discritize with finite elements and describe the solution procedure.
	
	\subsubsection{Augmented Lagrangian}
	
    The adjoint damage evolution for $b$ in \eqref{eq:adja} is challenging to efficiently solve. While the equation itself is linear, the $\dot{a}$ dependence makes the discretized operator dependent on the time-step. Therefore, we look to apply an augmented Lagrangian to cast this as a constant-operator global problem and a time-step dependent local problem. We introduce the auxillary field $z \in \mathcal{A}$, and constrain $z = \dot{a}{b}$ weakly through the Lagrange multiplier field $\chi \in L^2(\Omega)$. By writing the adjoint damage update as a minimization problem, we apply another augmented Lagrangian through the penalty parameter $r$ (See \ref{ap:adj}). This gives the adjoint damage evolution as 
    \begin{subequations}
    \begin{align}
    &0 = \int_{\Omega} \left[ \left( r(z - \dot{a} b) + \chi \right) \delta z + \frac{G_c \ell}{2 c_w} \nabla z \cdot \nabla \delta z \right] \, d \Omega && \forall \delta z \in \mathcal{A}, \label{eq:adjadmm1}\\
    & \dv{}{t} \left[ b D_a + \bar{\psi}^{* \prime \prime}  \dot{a} b  \right ] = \pdv{o}{a} + \pdv{^2 W^e}{a \partial \varepsilon} \cdot \nabla \xi + \dot{a} b \left( \pdv{^2 W^e}{a^2} + \frac{G_c }{4 \ell c_w} \pdv{^2 w^a}{a^2} \right)   \notag \\ 
    & \hspace{3cm} + \dot{a} b d^{\prime \prime}   \left[  W^p + \int_0^t g^* d\tau \right] - r(z - \dot{a} b) - \chi \qquad  &&\text{ on } \Omega, \label{eq:adjadmm3} \\
    &0 = \int_{\Omega} \left(  z -\dot{a} b \right) \delta \chi \, d \Omega &&  \forall \delta \chi \in L^2(\Omega). \label{eq:adjadmm2} 
    \end{align}
    \end{subequations}
    The first line \eqref{eq:adjadmm1} is linear constant-operator global problem for $z$. \eqref{eq:adjadmm3} is a linear local problem for $b$. Finally, the last line  \eqref{eq:adjadmm2} is the constraint that $z = \dot{a} b$ weakly.  We discuss the iterative method of solving this in the next section.
    
    \subsubsection{Discretization and Solution Procedure}
    
    The adjoint variables are discretized in the same manner as their forward counterparts. The adjoint displacement field $\xi$, the adjoint damage field $z$, and adjoint Lagrange multiplier fields are then
    \begin{equation}
    \xi = \sum_{i = 1}^{n_u} \xi_i N^u_i(x), \qquad z = \sum_{i = 1}^{n_a} z_i N^a_i(x), \qquad \chi =  \sum_{i = 1}^{n_a} \chi_i N^a_i(x).
    \end{equation} 
    The fields $b$, $\gamma$, and $\mu$ are discretized at quadrature points:
    \begin{equation}
    b(x_g) = b_g, \qquad \gamma(x_g) = \gamma_g, \qquad \mu (x_g) = \mu_g, 
    \end{equation}
    for some Gauss point $x_g$.
    The adjoint problem must be solved backwards in time. Similar to the forward problem, we use an explicit central difference scheme for the adjoint displacement variable. Then, we implicitly update the adjoint damage variables through an alternating direction method of multipliers. After these converge, the adjoint plastic variables are updated implicitly. For the $n + 1$ to the $n$ time-step the displacement updates are
    \begin{equation}
    \begin{aligned}
    \ddot{\xi}_i^{n + 1} &= M^{-1}_{ij} H_j^{n + 1}(u^{n + 1}, \varepsilon^{p,n+1}, \alpha^{n + 1}, \xi^{n + 1}, b^{n + 1}, \gamma^{n + 1}, \mu^{n + 1}),\\
    \dot{\xi}_i^{n + 1/2} &= \dot{\xi}_i^{n + 3/2} - \Delta t^{n + 1} \, \ddot{\xi}_i^{n + 1}, \\
    \xi_i^{n} &= \xi_i^{n + 1} - \Delta t^{n + 1/2} \, \dot{\xi}_i^{n + 1/2},
    \end{aligned}
    \end{equation}
    where,
    \begin{equation}
    H_j^n = \int_{\Omega} \left[ \left( - \nabla \xi^n \cdot \pdv{^2 W^e}{\varepsilon \partial \varepsilon}  - \dot{\alpha}^n b^n \pdv{^2 W^e}{\varepsilon \partial \alpha} - \gamma^n \dot{q}^n  \pdv{\bar{\sigma}_M}{\varepsilon} + \dot{q}^n \mu^n \cdot \pdv{M}{\varepsilon}\right) \cdot \nabla N^u_j  - \pdv{o}{u} \cdot N^u_j \right] \, d\Omega.
    \end{equation}
    The update for $b$ uses an implicit forward Euler method, coupled ADMM for fields $z$ and $\chi$. This results in iterations between a point-wise linear problem for $b$, a constant-matrix linear global problem for $z$, and an update for $\chi$. 
    
    We describe this for the $n + 1$ to $n$ time-step. Given $\xi^{n}$, intialize $\tilde{\chi}^0 = \chi^{n+1}$, $\tilde{z}^0 = z^{n+1}$, and iterate over~$i$:
    \begin{itemize}
    	\item {\textit{Step 1: Linear local problem.}} Update $\tilde{b}^{i+1}$ by solving at each $x_g$
    	\begin{equation}
    	\tilde{b}^{i + 1}_g = \frac{ \dot{\alpha}^{n+1}_g {b}^{n + 1}_g \left . \bar{\psi}^{* \prime \prime}  \right |_{t_{n + 1}} +  {b}^{n + 1}_g  \tilde{D}_{a,g}^{n+1} + \Delta t \left( r \tilde{z}^i(x_g) + \tilde{\chi}^i(x_g) - \left. \pdv{o}{a} \right |_{t_n} -   \left . \pdv{^2 W^e}{\alpha \partial \varepsilon} \right|_{t_n}  \cdot \nabla \xi^n \right)  }
    	{ \dot{\alpha}^n_g \left . \bar{\psi}^{* \prime \prime}  \right |_{t_n}  + \tilde{D}_{a,g}^n + \dot{\alpha}^n_g \left(  \left [ \pdv{^2 W^e}{\alpha^2} + \frac{G_c }{4 \ell c_w } \pdv{^2w^a}{a^2} \right]_{t_n} +  d^{\prime \prime}   \left[  W^p + \int_0^t g^* d\tau \right]_{t_n} + r\right) },
    	\end{equation}
    	where
    	\begin{equation}
    	\tilde{D}_{a,g}^n = \left [ \pdv{W^e}{\alpha}  + \frac{G_c }{4 \ell c_w }\pdv{w^a}{\alpha} + \pdv{d}{\alpha} \left(W^p + \int_0^t g^*\, d\tau \right) + \pdv{\bar{\psi}^*}{\dot{\alpha}}  \right ]_{x_g, t_n} - r(a^n|_{x_g} - \alpha^n_g) - \lambda^n |_{x_g}.
    	\end{equation}
    	
    	\item {\textit{Step 2: Linear global problem.}} Update $\tilde{z}^{i+1}$ by solving
    	\begin{equation}
    	{K}_{pj} \, \tilde{z}^{i+1}_j = U_p(\tilde{b}^{i+1}, \tilde{\chi}^i),
    	\end{equation}
    	where
    	\begin{equation}
    	U_p(b, \chi) = \int_{\Omega} \left( r \dot{\alpha}^n b - \chi \right) N^a_p \, d\Omega.
    	\end{equation}
    	
    	\item {\textit{Step 3: Update Lagrange multiplier.}} Update $\tilde{\chi}^{i+1}$ by
    	\begin{equation}
    	\tilde{\chi}^{i + 1}_j = \tilde{\chi}^i_j + r(\tilde{z}^{i + 1}_j - S^{-1}_{jk} \hat{z}^{i + 1}_k ),
    	\end{equation}
    	where
    	\begin{equation}
    	\hat{z}^{i+1}_k = \int_{\Omega} \dot{\alpha}^n \tilde{b}^{i+1} N^a_k \, d\Omega.
    	\end{equation}
    	Note: this is the weak form of the update $\Delta \chi = r(z - \dot{\alpha} b )$.
    	
    	\item {\textit{Step 4: Check for convergence.}} Check both primal and dual feasibility,
    	\begin{equation}
    	\begin{aligned}
    	r_p &:= \norm{ \bar{z}^{i+1} - \hat{z}^{i+1}}_{l^2} \leq \frac{1}{\sqrt{n_a}} r^{tol}_{abs} + r^{tol}_{rel} \max \left(\norm{\hat{z}^{i+1}}_{l^2}, \norm{\bar{z}^{i+1}}_{l^2}\right), \\
    	r_d &:= r \norm{\bar{z}^{i+1} - \bar{z}^{i}}_{l^2} \leq \frac{1}{\sqrt{n_a}} r^{tol}_{abs} + r^{tol}_{rel} \norm{\bar{\chi}^{i + 1}},
    	\end{aligned}
    	\end{equation} 
    	where
    	\begin{equation}
    	\bar{z}^{i + 1}_j = S_{jk} \tilde{z}^{i +1}_k, \qquad  \bar{\chi}^{i + 1}_j = S_{jk} \tilde{\chi}^{i +1}_k.
    	\end{equation}
    \end{itemize}
    until convergence, and set $b^{n} = \tilde{b}^i$, $z^{n} = \tilde{z}^i$, and $\chi^{n} = \tilde{\chi}^i$. We adapt the penalty value $r$ similarly to the forward problem in \eqref{eq:updater}.
    
    Finally, the adjoint plastic variables $\gamma^{n}_g$ and $\mu^{n}_g$ are implicity updated by solving at each quadrature point: 
    \begin{equation}
    \begin{aligned}
    & \left[ \gamma \left( \bar{\sigma}_M - \sigma_0 - \pdv{\bar{\psi}^*}{\dot{q}}\right)- \gamma \dot{q} \pdv{^2 \bar{g}^*}{\dot{q}^2} + \pdv{\bar{g}^*}{\dot{q}} \left( \int_t^T b \dot{\alpha} d^{\prime}(\alpha) d\tau\right) - \mu \cdot M \right]_{\substack{t = t_{n + 1} \\ x = x_g \ \ \ }} -  \\
    &\qquad \qquad  \left[ \gamma \left( \bar{\sigma}_M - \sigma_0 - \pdv{\bar{\psi}^*}{\dot{q}}\right)- \gamma \dot{q} \pdv{^2 \bar{g}^*}{\dot{q}^2} + \pdv{\bar{g}^*}{\dot{q}} \left( \int_t^T b \dot{\alpha} d^{\prime}(\alpha) d\tau\right) - \mu \cdot M \right]_{\substack{t = t_n \\ x = x_g}}  \\
    &\qquad \qquad \    = \Delta t \left( \left. \pdv{o}{q} \right |_{\substack{t = t_n \\ x = x_g}} + b_g^{n} \dot{\alpha}_g^n d^{\prime}(\alpha_g^n) \left. \pdv{W^p}{q} \right |_{\substack{t = t_n \\ x = x_g}} - \gamma_g^n \dot{q}_g^n  \left. \pdv{\sigma_0}{q} \right|_{\substack{t = t_n \\ x = x_g}} \right ) \\
    & \mu_g^{n + 1} - \mu_g^n = \Delta t \left( \left. \pdv{o}{\varepsilon^p} \right |_{\substack{t = t_n \\ x = x_g}} + \nabla \xi^n \cdot \left. \pdv{^2 W^e}{\varepsilon \partial \varepsilon^p} \right |_{\substack{t = t_n \\ x = x_g}} \right. \\
    & \qquad \qquad \qquad \qquad  \left. + b^n_g \dot{\alpha}^n_g \left. \pdv{^2 W^e}{ \alpha \partial \varepsilon^p} \right |_{\substack{t = t_n \\ x = x_g}}+ \gamma \dot{q}^n_g \left. \pdv{\bar{\sigma}_M}{\varepsilon^p} \right |_{\substack{t = t_n \\ x = x_g}} - \dot{q}_g^n \mu^n_g \cdot \left. \pdv{M}{\varepsilon^p} \right |_{\substack{t = t_n \\ x = x_g}} \right).  \\        
    \end{aligned}
    \end{equation}
    This is a linear system of equations which may be solved by direct inversion.
    
	\subsection{Sensitivities and Design Updates}
\begin{figure}
		\begin{center}
				\includegraphics[width=0.95\textwidth]{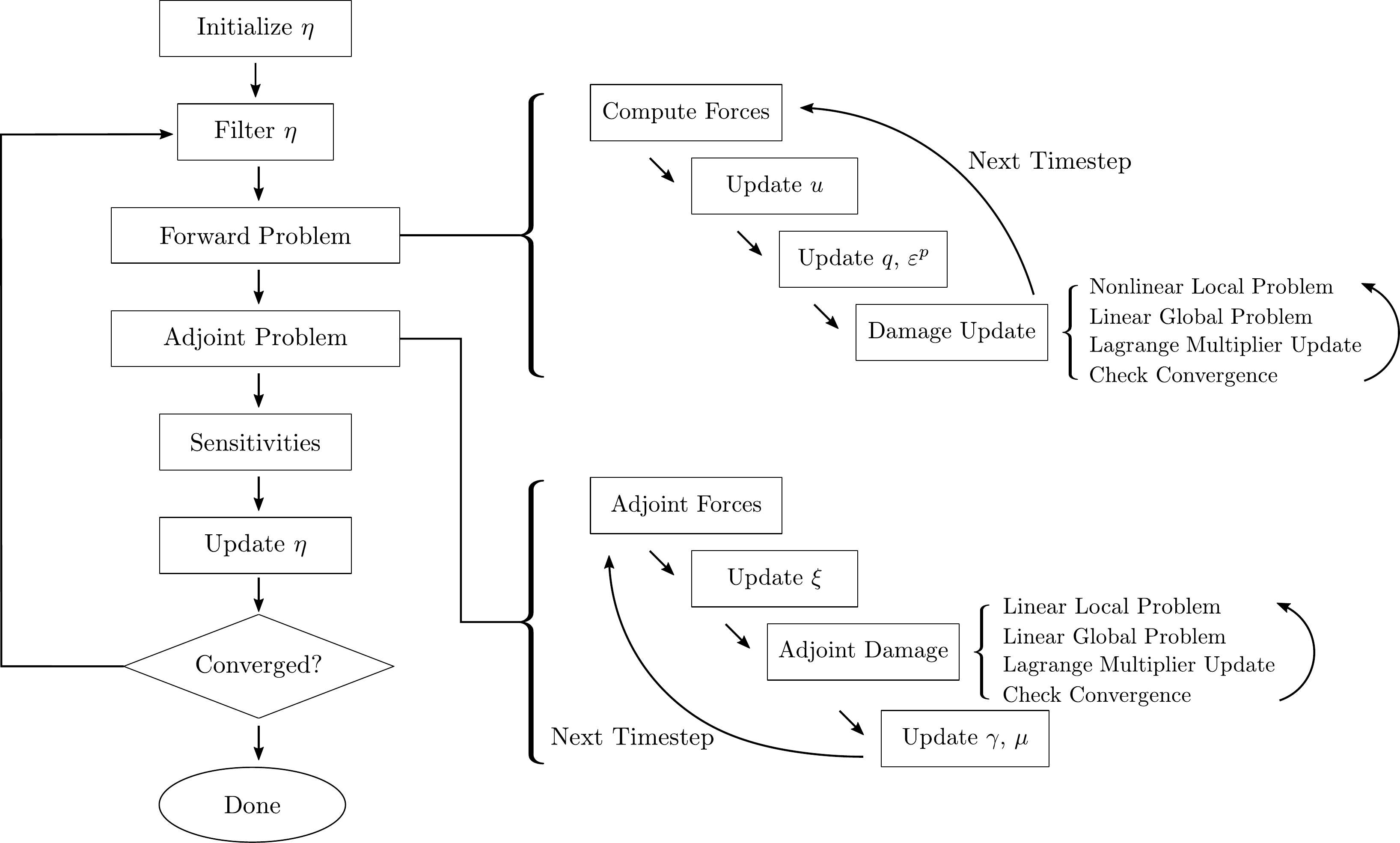}
		\end{center}
		\caption{Diagram of the computational method for gradient-based topology optimization over the dynamic trajectory with plasticity and damage. }
		\label{fig:diagram}
	\end{figure}

	Optimal design problems in structural mechanics often lead to ill-posed minimization problems, where minimizing sequences develop fine scale oscillations~\cite{Kohn1986,strang1986}. To recover a well-posed problem, we filter the design variable $\eta$. These density-based filtering methods have been shown to lead to well-posed problems for linear, static compliance optimization. We consider $\eta$ constant on each element, and adopt a discrete re-normalized filter with a linear weight function~\cite{Bourdin2001}. Sensitivities, accounting for the filtering, are then computed from \eqref{eq:sensitivities}. These are used to update $\eta$ using the gradient-based method of moving asymptotes~(MMA)~\cite{K.1987}. This process is continued until convergence. Figure~\ref{fig:diagram} shows a flow diagram of the entire computational process. 
	
	\section{Material Interpolation} \label{sec:interp}
	In the preceding section, we developed a computational method for evolving the forward and adjoint problem to compute sensitivities. However, we still must define how the material parameters depend on the design parameter $\eta$. That is, we must determine how the material density, elastic energy, plastic potential and dissipation, and also the damage parameters depend on $\eta$. In this section, we discuss interpolation schemes for both solid-void designs, as well as designs composed of two materials of differing parameters. 
	
	\subsection{Solid-Void Designs}
	We consider $\eta$ as differentiating between void at $\eta = \eta_{min} << 1$ and solid at $\eta = 1$. Similar to traditional topology optimization, we would like to penalize intermediate densities so converged designs are dominated by regions of completely solid or void. In the following, the subscript $0$ denotes parameters for the completely solid material. We propose the following interpolation scheme:
	\paragraph{Material density} We consider $\eta$ a density variable, and assume the material density varies linearly: 
    \begin{equation}
        \rho(\eta) = \eta \rho_0.
    \end{equation}
    \paragraph{Elastic Energy} For simplicity, we consider a separable dependence for the elastic energy through a Bezier curve interpolation. This ensures that the ratio of stiffness to density does not go to zero in the limit of small $\eta$. This mitigates spurious dynamical modes which could arise from artificial acoustic properties of the voids~\cite{Bendsoee2004}. We consider
    \begin{equation}
        W^e(\varepsilon, \varepsilon^p, a, \eta) = B_e(\eta) W^e_0(\varepsilon, \varepsilon^p, a),
    \end{equation}
    where $W^e_0$ is the elastic energy of the solid and $B_e(\eta)$ is defined through
    \begin{equation}
    \begin{aligned}
        \eta &= \frac{1 - k_2}{k_1 - k_2} (3v - 3v^2) + v^3, \\
        B_e &= k_1 \frac{1 -k_2}{k_1 - k_2} (3v - 3v^2) + v^3.
    \end{aligned}
    \end{equation}
    Given $\eta$, the top equation may be solved for $v$, which is then used to compute $B_e$ in the second equation. $k_1$ and $k_2$ are the derivative values $\frac{dB}{d \eta}$ at $\eta = 0$ and $\eta = 1$, respectively. 
    \paragraph{Plastic potentials}For the plastic potentials, we will again consider a separable dependence 
    \begin{equation}
        W^p(q, \eta) = B_p(\eta) W^p_0(q), \qquad g^*(\dot{q}, \eta) = B_p(\eta) g^*_0(\dot{q}).
    \end{equation}
    However, care must be taken in choosing $B_p(\eta)$, as we require this interpolation to satisfy certain properties:
    \begin{itemize}
        \item \textit{Strong voids} : The yield stress should be sufficiently high as to reduce excessive permanent deformation in the void regions. Additionally, we do not want to waste computational effort on plastic updates in the voids. This requires
        \begin{equation}
            1 < \frac{B_p(\eta_{min})}{B_e(\eta_{min})}.
        \end{equation}
        \item \textit{Unfavorable intermediate densities} : The interpolation of the plastic potential should ensure that the relative yield stress is not excessively high in regions of intermediate density, so optimal solutions are dominated by regions of either completely solid or void. This requires
        \begin{equation}
            \frac{B_p(\eta)}{B_e(\eta)} < \tau_p \quad \forall \eta \in [\eta_1, \eta_2],
        \end{equation}
        where $\eta_{min} < \eta_1 < \eta_2 < 1$ and $\tau_p \sim 1$.
    \end{itemize}
    We may accomplish both of these by considering a shifted Bezier curve interpolation as
    \begin{equation}
        B_p(\eta) = \frac{B_e(\eta) + \delta_p }{1 + \delta_p},
    \end{equation}
    where $B_e(\eta_{min}) < \delta_p << 1$. 
    \paragraph{Damage parameters} We now discuss the interpolation for the damage behavior. For simplicity, we assume that the normalized damage potential $w^a(a)$ is independent of the density. The damage length scale will also be considered constant with density
    \begin{equation}
        \ell(\eta) = \ell_0.
    \end{equation}
    This allows the same computational mesh to resolve damage in both the solid and void regions. Then, we must only prescribe the interpolation on the toughness $G_c$. We assume a separable dependency
    \begin{equation}
        G_c(\eta) = B_a(\eta) G_{c0},
    \end{equation}
    where the interpolation function $B_a$ must satisfy the following:
    \begin{itemize}
        \item \textit{Boundary condition preservation} : The behavior at the solid-void interface should be nearly equivalent to the natural boundary conditions. This ensures that the voids behave similarly to free boundaries and do not add artificial toughness. This requires 
        \begin{equation}
            B_a(\eta_{min}) << B_a(1)
        \end{equation}
        \item \textit{Tough voids}: We require that the damage not propagate through the void regions, which could result in damage ''jumping" from one solid region to another by moving through voids. This requires 
        \begin{equation}
            1 < \frac{B_p(\eta_{min})}{B_e(\eta_{min})} << \frac{B_a(\eta_{min})}{B_e(\eta_{min})},
        \end{equation}
        ensuring that the relative toughness of the voids is much larger than that of the solid.
        \item \textit{Unfavorable intermediate densities} : The damage interpolation should ensure that the relative toughness is not excessively high in regions of intermediate density, so optimal solutions are dominated by regions of either completely solid or void. This requires
        \begin{equation}
            \frac{B_a(\eta)}{B_e(\eta)} < \tau_a \quad \forall \eta \in [\eta_1, \eta_2],
        \end{equation}
        where $\eta_{min} < \eta_1 < \eta_2 < 1$ and $\tau_a \sim 1$.
    \end{itemize}
    We may again accomplish these through a shifted Bezier curve,
    \begin{equation}
        B_a(\eta) = \frac{B_e(\eta) + \delta_a }{1 + \delta_a},
    \end{equation}
    where $B_e(\eta_{min}) << \delta_p < \delta_a << 1$.
    
\begin{figure}
		\begin{center}
			\includegraphics[width=0.6\textwidth]{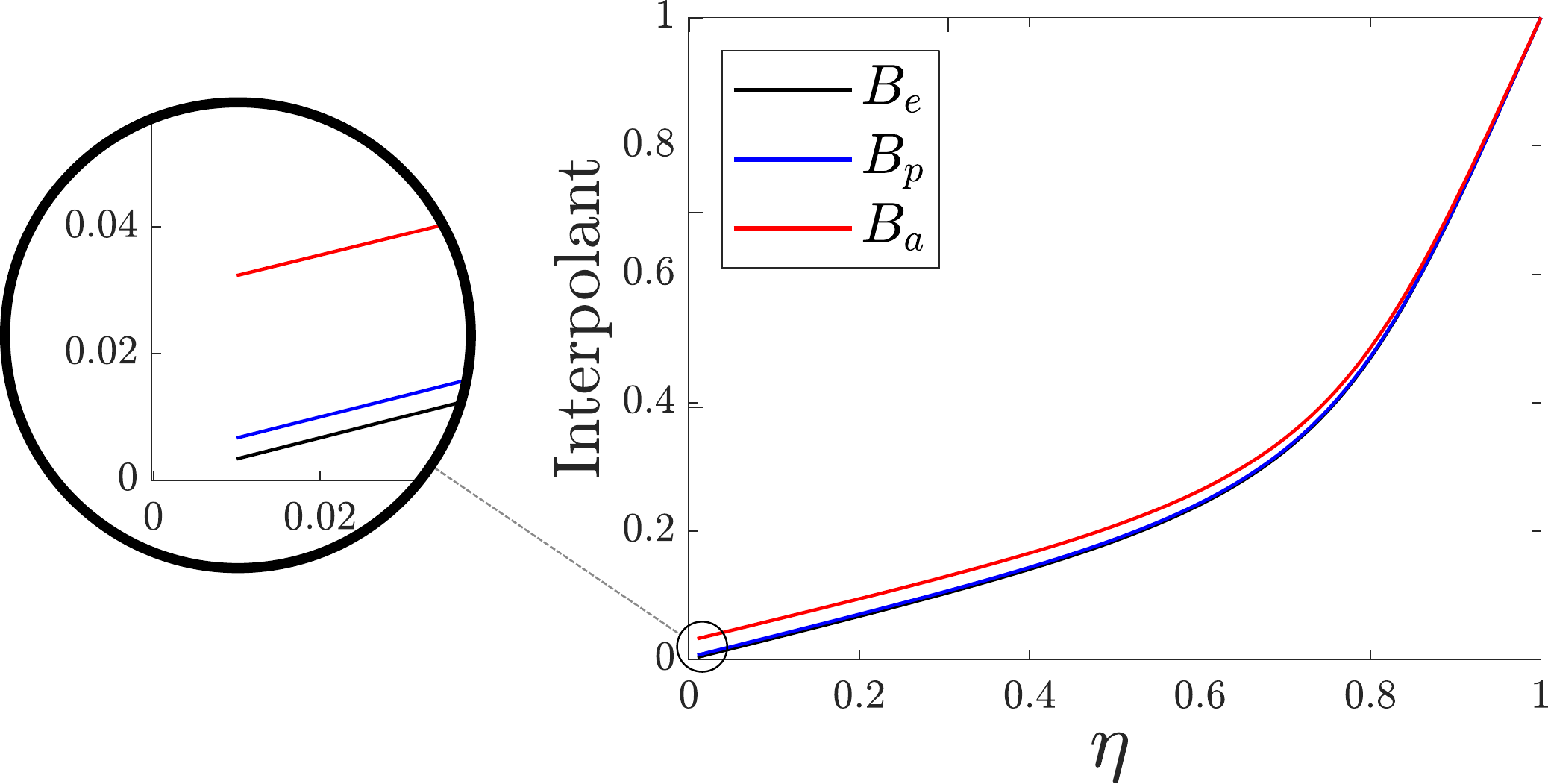}
		\end{center}
		\caption{Plot of the interpolation functions for the elasticity ($B_e$), plasticity ($B_p$), and damage ($B_a$) for parameters $k_1 = 0.2$, $k_2 = 5.0$, $\eta_{min} = 0.01$, $\delta_p = k_1 \eta_{min}$, and $\delta_a = 9 k_1 \eta_{min}$. Here, $\eta_{min} = 0.01$.}
		\label{fig:interpo}
	\end{figure}

    For our investigation, we choose a value of $\delta_p = k_1 \eta_{min}$, $\delta_a = 9 k_1 \eta_{min}$. Thus, the yield strain of the void regions is roughly twice that of the solid. Additionally, the voids have around $10$ times the relative toughness of the solid regions. Figure~\ref{fig:interpo} shows these interpolation functions plotted for typical values. 
    
    \subsection{Two-Material Design}
	We now consider designs composed of two materials, where $\eta = 0$ and $\eta = 1$  represents solids of either species. We propose to interpolate the majority of these parameters through standard power-law functions. These penalize regions of intermediate densities, while also being efficient and simple to implement. 
	
	\section{Examples} \label{sec:results}
	
	We now demonstrate the methodology using two examples.  The examples are also of independent interest for the insights they offer on damage resistant structures. The first is a solid-void design to resist impulse loading. The second example explores the trade-offs between strength and toughness in a spall-resistant structure composed of two different materials undergoing dynamic impact.   The forward dynamics, adjoint problem, sensitivity calculation, and MMA update schemes are implemented using the deal.II C++ finite element library~\cite{Bangerth2007}.

	\subsection{\bfseries Solid-Void for Blast Loading}
\begin{figure}
		\begin{center}
			\includegraphics[width=0.6\textwidth]{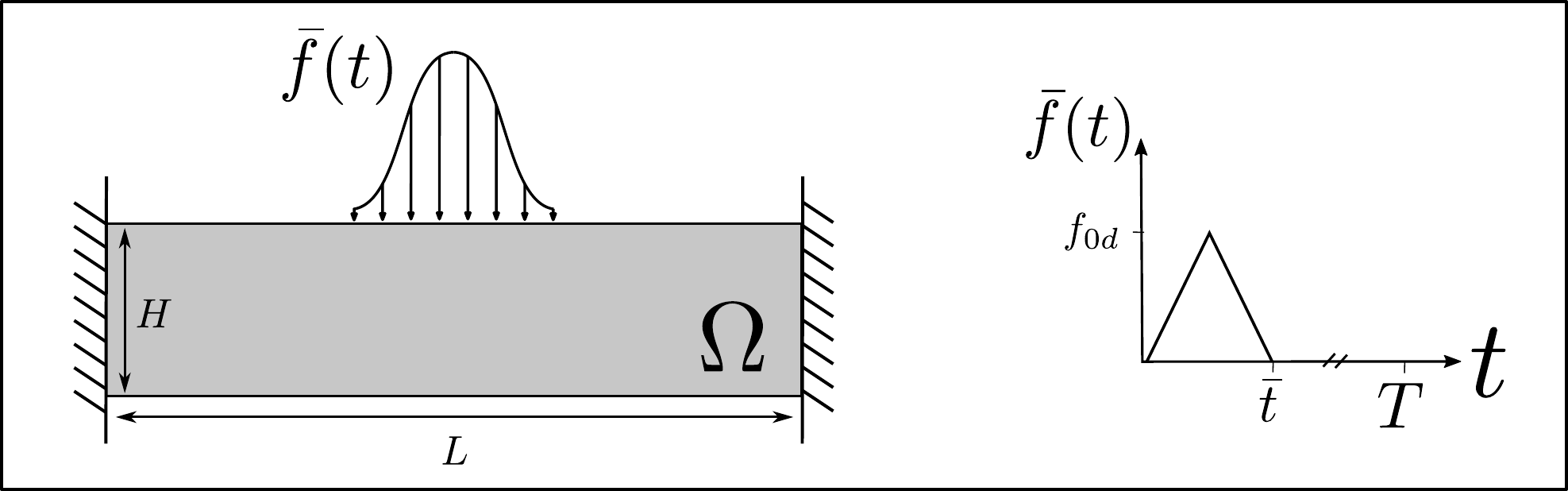}
		\end{center}
		\caption{Geometry and dynamic impulse loading we consider for the solid-void structure.}
		\label{fig:solidVoid}
	\end{figure}

	We consider $\eta$ as a density variable distinguishing between solid material and void. To model blast loading, we assume a fixed loading prescribed on the boundary. Thus, we consider a rectangular 2D geometry and impulse loading as shown in Figure~\ref{fig:solidVoid}. We look to minimize a sum of the time-space norm of the displacements, plastic dissipation, and damage dissipation
	\begin{equation} \label{eq:objective}
		\mathcal{O} = \frac{\sigma_{y0} L }{T^{1/s}}  \| \left( \| u \|_{H^1(\Omega)} \right) \|_{L^s(0, T)}  + c_p \mathcal{D}^p + c_a \mathcal{D}^a,
	\end{equation}
	where $c_p$ and $c_a$ are weights, and $\mathcal{D}^p$ and $\mathcal{D}^a$ are measures of the dissipated energy to plasticity and damage, 
	\begin{equation}
	    \mathcal{D}^p = \int_\Omega d(a(T)) \left( \left. \tilde{W}^p \right |_{t = T}+ \int_{0}^T \tilde{g}^* \ dt \right) \ d\Omega, \qquad \mathcal{D}^a = \int_\Omega \left [  \frac{\tilde{G}_c w^a(a(T))}{4 c_w \ell} + \int_0^T \tilde{\psi}^* \ dt \right]  \ d\Omega.
	\end{equation}
	We use a modified interpolation scheme in the objective to penalize intermediate densities. That is, we choose $\tilde{W}^p$, $\tilde{g}^*$,  $\tilde{G}_c$, and $\tilde{\psi}^*$ to remain relatively large for intermediate $\eta$. Thus, we consider a concave power-law interpolation
	\begin{equation}
	    \tilde{W}^p(q, \eta) = P(\eta) W_0^p(q), \quad g^*(\dot{q}, \eta) = P(\eta) g_0^*(q), \quad \tilde{G}_c(\eta) = P(\eta) G_{c0}, \quad   \psi^*(\dot{q}, \eta) = P(\eta) \psi_0^*(q),
	\end{equation}
	where
	\begin{equation}
	    P(\eta) = 1 - (1 - \eta)^{p_O}.
	\end{equation}
	Here, $p_O$ is a growth factor parameter. $s$ is the power for the norm in time. Because we intend minimize the largest displacements, we choose $s = 4$ for the following studies.
	
	We consider material parameters shown in Table~\ref{tab:solidvoidparams}. We consider a Gaussian loading profile of standard deviation $L/20$, truncated to a total width of $L/5$. We use objective penalty values of $c_p = 5, \ c_a = 50$. Thus, we look to heavily penalize damage. For the interpolation parameters, we linearly update the Bezier slopes from $k_1 = 0.5$, $k_2 = 2.0$ to $k_1 = 0.125$, $k_2 = 8$ from the first to the $50$th iteration. This allows the structure topology to more free change at lower iterations before intermediate densities are severely penalized, and is standard practice in topology optimization~\cite{Bendsoee2004}. Because the structure may not be able to withstand the loading (without severe damage) for the early iterations, we begin with a lower loading amplitude before gradually increasing to the final desired value. We set the loading amplitude to be $70\%$ of the final value until iteration $60$, which we then linearly increase to the final value by iteration $100$. Computations are performed on a $100 \times 25$ mesh, with a density filter radius of $0.021 L$. Additionally, we restrict the amount of material used to be no more than half the volume of $\Omega$. Designs are then initialized to uniform density fields equal to the total allowed volume fraction $\eta = 0.5$. We consider designs converged when the maximum change in density variable is less than $10^{-3}$, or after $300$ iterations.
	
	\begin{table}
    \begin{center}
    \begin{tabular}{ |p{3cm}||p{3cm}|p{6cm}| }
    \hline
    Parameter & Value & Description\\
    \hline
    \multicolumn{3}{|c|}{\textbf{Geometric Parameters}} \\
    \hline
    $L$ & $1$  & Length of domain \\
    $H$ & $0.25$  &  Height of domain \\
    \hline
    
    \multicolumn{3}{|c|}{\textbf{Elastic Material Parameters}} \\
    \hline
    $E$ & $0.5$  & Young's modulus \\
    $\nu$ & $0.3$  & Poisson ratio \\
    $\rho$ & 0.05 & Density \\
    \hline
    \multicolumn{3}{|c|}{\textbf{Plastic Material Parameters}} \\
    \hline
    $\sigma_{y0}$ & $0.01 E$  & Yield strength \\
    $\varepsilon^p_0$ & $0.1$  & Reference plastic strain \\
    $n$  & $10$ & Isotropic hardening power \\
    $\dot{\varepsilon}^p_0$ & $1.0$  & Reference plastic strain rate\\
    $m$  & $6$ & Rate sensitivity power\\
    \hline
    \multicolumn{3}{|c|}{\textbf{Damage Material Parameters}} \\
    \hline
    $G_{c0}$ & $1.5\times 10^{-4}$  & Toughness \\
    $\ell$ & $0.02$  & Damage length scale \\
    $d_1$ & $0.01$ & Relative stiffness when fully damaged  \\
    $w_1$ & $0.95$  & Damage hardening parameter \\
    \hline
    \end{tabular}\\
    \caption{Geometric and material parameters for the solid-void structures.}
    \label{tab:solidvoidparams}
    \end{center}
    \end{table}
	
\begin{figure}
		\begin{center}
			\includegraphics[width=0.95\textwidth]{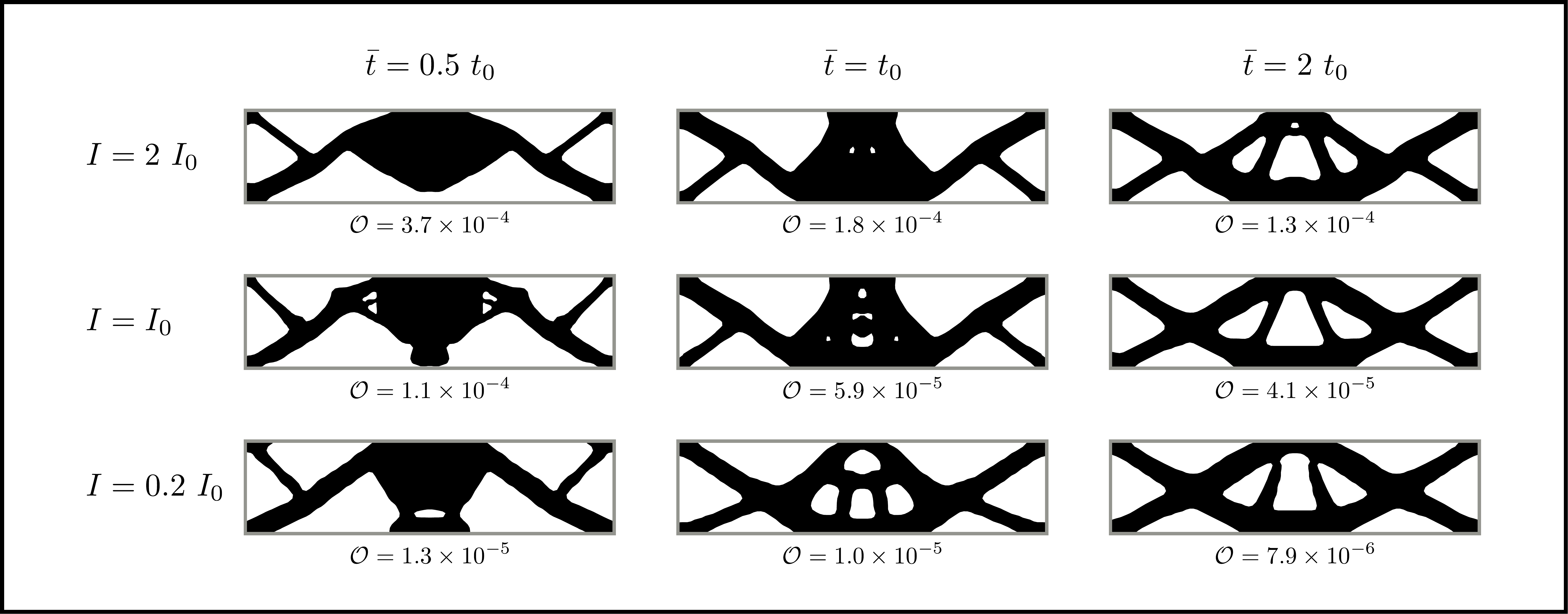}
		\end{center}
		\caption{Converged solid-void designs under impulse loading following contour smoothing. Along each row, the impulse is constant, while along the columns we vary the loading duration. Values of the objective are shown for each of the designs. All of the designs saturated the constraint that $V \leq 0.5 | \Omega |$. }
		\label{fig:solid_void_designs}
	\end{figure}

	We explore optimal designs for varying impulse magnitude and loading duration. Here, we consider a reference loading duration $t_0 = 1.47 \ L/c_L$, roughly the time that it takes a longitudinal wave to traverse three half-length of the domain. For the reference impulse, we consider $I_0 = 6.3\times10^{-4}\ L^2\sqrt{E \rho}$. The simulation time is set to $T = 19 t_0$. Figure~\ref{fig:solid_void_designs} shows the converged design after contour smoothing in MATLAB\textsuperscript{\textregistered}. Along each row, the loading impulse is constant, while along each columns the loading duration is constant. Although the structures share similar supports near the boundaries, their topologies near the loading site vary drastically. We see that for that for longer loading duration (right column), the structure is similar to what we would expect from static compliance optimization: truss-like members forming triangular structures~\cite{Bendsoee2004}. However, for shorter loading duration, the structures have more mass congregated underneath the applied load. This not only provides damage resistance, but the additional inertia also reduces the energy the structure absorbs from the impulse loading. We also see more mass placed near the loading surface for large impulse magnitude. This is likely to reduce plasticity and damage near the loading site.

	
	\subsection{\bfseries Two Material Design for Impact}
\begin{figure}
		\begin{center}
			\includegraphics[width=0.95\textwidth]{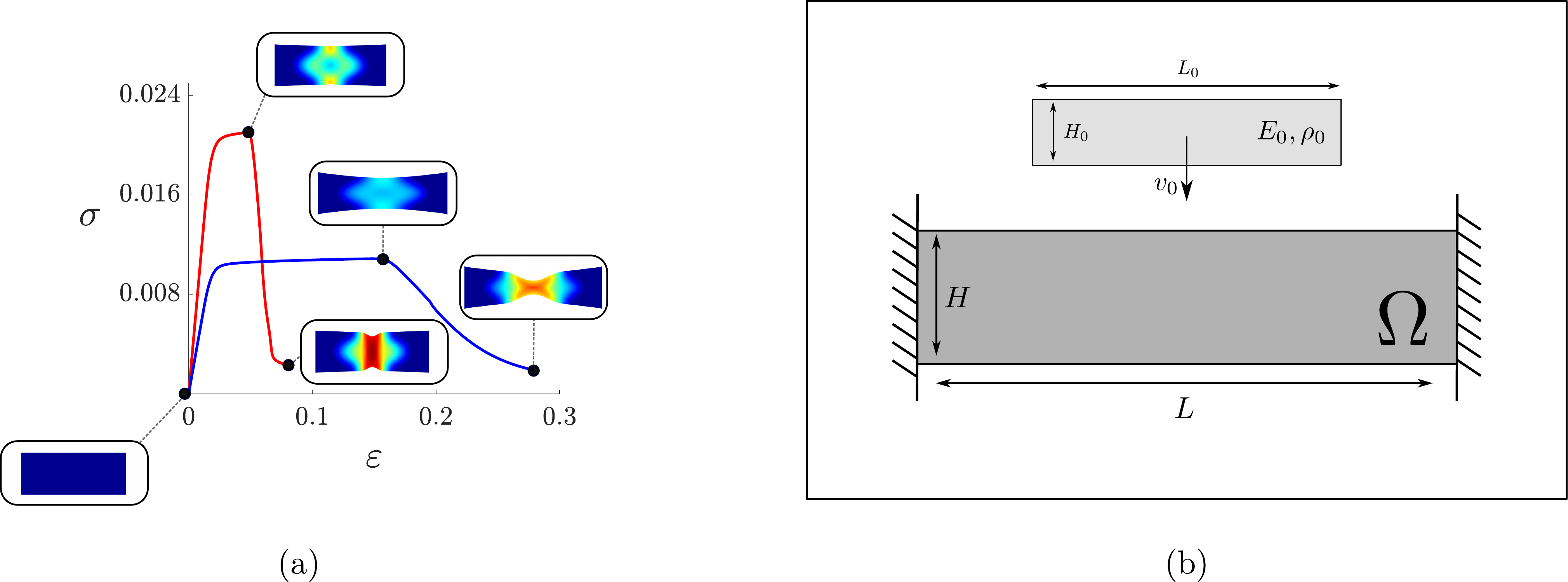}
		\end{center}
		\caption{(a) Normalized stress-strain response of the strong (red) and tough (blue) material in a uniaxial quasi-static 2D tensile test. Damage fields are plotted on deformed configurations at a few points throughout loading. (b) Geometry and loading for the two-material structure. }
		\label{fig:strongTough}
	\end{figure}

	We now consider the design of a structure composed of two materials undergoing impact. Figure~\ref{fig:strongTough}\hyperref[fig:strongTough]{a} shows the stress-strain response of the two materials in a quasi-static tensile simulation. While one material has double the strength and stiffness (red curve), the other has roughly twice the toughness (blue curve). We represent the strong solid with  $\eta = 1$, and the tough solid with $\eta = 0$. Thus, we consider,
    \begin{equation}
        E_1 < E_2, \quad \left( \sigma_{y0} \right)_1 < \left( \sigma_{y0} \right)_2, \quad \left(G_c \right)_1 > \left( G_c \right)_2,
    \end{equation}
    where $E$, $G_c$, and $\sigma_{y0}$ denotes the elastic modulus, fracture toughness, and yield stress. The subscripts $1$ and $2$ denotes properties of the tough and strong solid, respectively. For simplicity, we assume the rest of the material properties are identical (density, hardening parameters, damage length scale). As discussed in the previous section, we adopt a power-law interpolation for material parameters. However, to ensure that the intermediate $\eta$ remains unfavorable, we must carefully choose the concavity of each of the interpolation functions. Since it is assumed that a larger value for each of the differing parameters is favorable, the interpolation is convex for all of these:
    \begin{equation}
    \begin{aligned}
        E(\eta) &=  E_1 + \eta^p (E_2 - E_1), \\
        \sigma_{y0} (\eta) &=  \left( \sigma_{y0} \right)_1 + \eta^p \left[ \left( \sigma_{y0} \right)_2  - \left( \sigma_{y0} \right)_1 \right ]\\
        G_c(\eta) &= \left( G_c \right)_2 +  \left( 1 - \eta \right)^p \left[ \left( G_c \right)_1 - \left( G_c \right)_2 \right] .
    \end{aligned}
    \end{equation}
    Similarly to the solid-void structure, we start with a penalty value of $p = 2$, and linearly increase it to $p = 8$ by the $100$th iteration and onward. 
    
    We consider the geometry as Figure~\ref{fig:strongTough}\hyperref[fig:strongTough]{b}. Here, we consider a linear elastic flyer of density $\rho_0$ and elastic modulus $E_0$ with an initial velocity of $v_0$.    We note that enforcing strict contact conditions would complicate the adjoint sensitivity calculations, and also be computationally expense.  Therefore, we consider a relaxation by introducing a layer of asymmetric linear elastic elements between the domain $\Omega$ and the flyer.  These elements have a high bulk modulus  in compression, with nearly zero resistance to shear or hydrostatic tension. Therefore, they may support compressive contact forces, while allowing the flyer and substrate to separate. This is consistent with the adjoint formulation we have derived, while providing the necessary physics of contact and separation. However, we are limited to cases where the impact site is known \textit{a priori} and the impacting faces are parallel.  
    
\begin{figure}
		\begin{center}
			\includegraphics[width=0.9\textwidth]{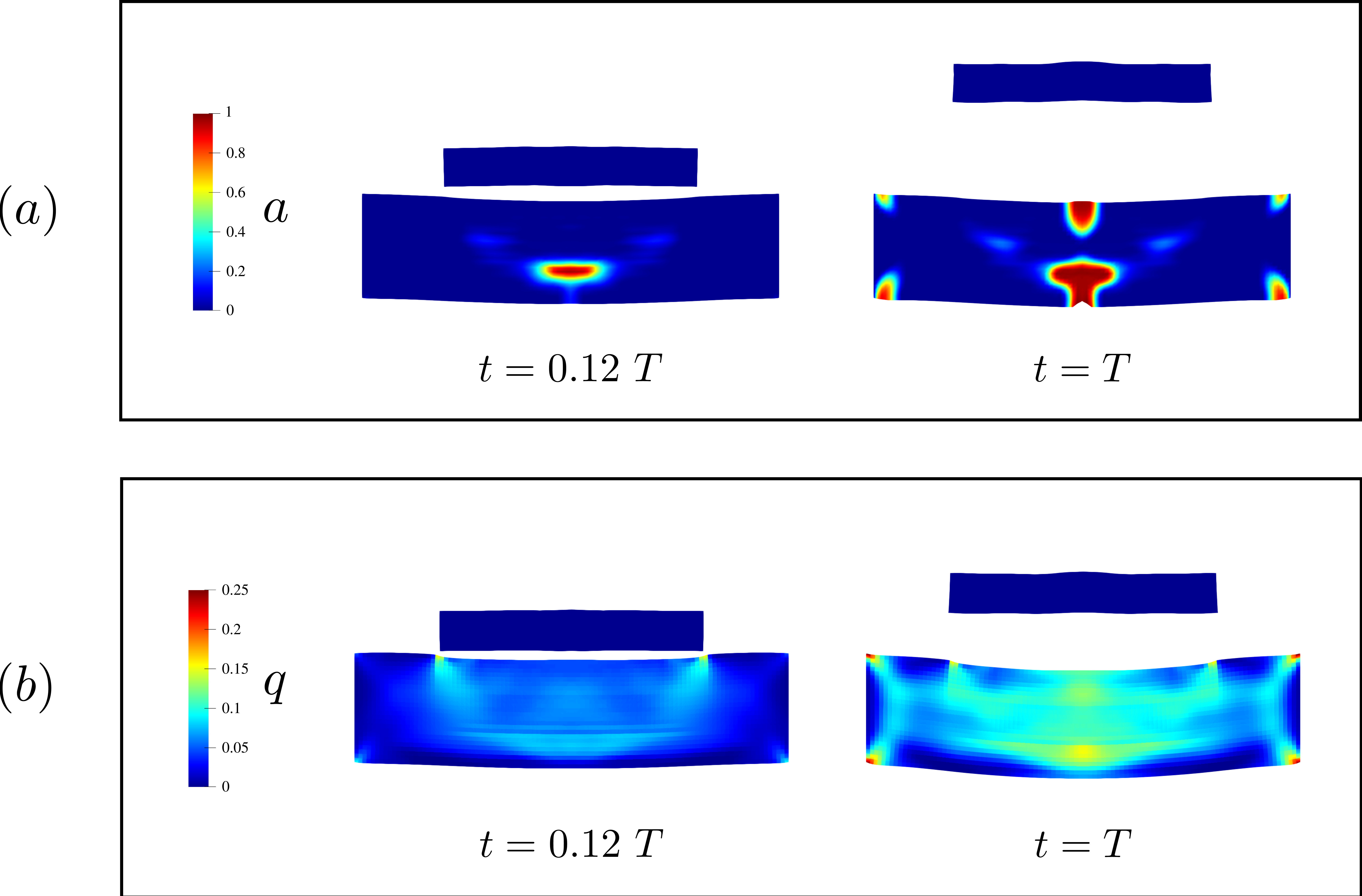}
		\end{center}
		\caption{(a) Damage field of a domain composed entirely of the strong material shortly after impact (left), and long after impact (right). (b) Accumulated plasticity field of a domain composed entirely of the tough material shortly after impact (left), and long after impact (right). }
		\label{fig:wholets}
	\end{figure}

	First, we consider the optimal design of the multi-material structure undergoing a relatively high impact velocity of $0.110  c_L$, where $c_L$ is the longitudinal wave speed of the strong material. Here, we consider a simulation time of $T = 6.5 \ L/c_L$. In this case, a structure composed entirely of strong material experiences heavy damage, as shown in Figure~\ref{fig:wholets}\hyperref[fig:wholets]{a}. The damage nucleates internally along a line parallel to the impact surface, which is characteristic of spall failure. Conversely, a structure made of only the tough material has large permanent deformation. There is regions of plasticity near the impact site as well as hinging near the boundary, as shown in Figure~\ref{fig:wholets}\hyperref[fig:wholets]{b}. We apply the optimal design approach to this loading scenario, as we hypothesize that a mixture of both strong and tough material will yield a structure of better performance. We consider the objective shown in \eqref{eq:objective}. Since $\eta$ does not have a large effect on the dissipation functions, we do not need to modify the interpolation scheme in the objective as was done in the solid-void case. Table~\ref{tab:multiparams} shows the material parameters we consider. Computations are performed on a $100\times25$ mesh for the rectangular domain, with a $60\times16$ mesh for the flyer. A filter radius of $0.021L$ is used. Again, we use objective weights of $c_p = 5$, $c_a = 50$. The design is initialized to a uniform density field of $\eta = 0.5$. The converged optimal design is shown in Figure~\ref{fig:strongTough}\hyperref[fig:strongTough]{a}. Red regions are occupied by strong material, and blue regions by the tough material. We see regions of strong (red) material near the boundaries and the impact site to mitigate large deformations and provide strength. However, the center is occupied by tough (blue) material to control spall. In terms of quantified performance, the converged design yields an objective value of $\mathcal{O} = 10.6 \times 10^{-4}$. This is improved performance over both the completely strong structure ($\mathcal{O} = 29.7 \times 10^{-4}$), and the completely tough structure ($\mathcal{O} = 11.3 \times 10^{-4}$).
	
	\begin{table}
    \begin{center}
    \begin{tabular}{ |p{3cm}||p{3cm}|p{6cm}| }
    \hline
    Parameter & Value & Description\\
    \hline
    \multicolumn{3}{|c|}{\textbf{Geometric Parameters}} \\
    \hline
    $L$ & $1$  & Length of domain \\
    $H$ & $0.25$  &  Height of domain \\
    $L_0$ & $0.6$  & Length of flyer \\
    $H_0$ & $0.1$  &  Height of flyer \\
    \hline
    
    \multicolumn{3}{|c|}{\textbf{Elastic Material Parameters}} \\
    \hline
    $E_1$ & $0.5$  & Young's modulus of tough material\\
    $E_2$ & $1.0$  & Young's modulus of strong material\\
    $\nu$ & $0.3$  & Poisson ratio \\
    $\rho$ & 0.05 & Density\\
    $E_0$  & 0.3 & Young's modulus of flyer \\
    $\nu_0$ & 0.4 & Poisson ratio of flyer \\
    $\rho_0$ & 0.02 & Denstiy of flyer \\
    \hline
    \multicolumn{3}{|c|}{\textbf{Plastic Material Parameters}} \\
    \hline
    $(\sigma_{y0})_1$ & $0.5\times10^{-2}$  & Yield strength of tough material\\
    $(\sigma_{y0})_2$ & $1.0\times10^{-2}$  & Yield strength of strong material\\
    $\varepsilon^p_0$ & $0.1$  & Reference plastic strain \\
    $n$  & $3$ & Isotropic hardening power \\
    $\dot{\varepsilon}^p_0$ & $1.0$  & Reference plastic strain rate\\
    $m$  & $3$ & Rate sensitivity power\\
    \hline
    \multicolumn{3}{|c|}{\textbf{Damage Material Parameters}} \\
    \hline
    $(G_{c0})_1$ & $1.0\times 10^{-4}$  & Toughness of tough material\\
    $(G_{c0})_2$ & $0.5\times 10^{-4}$  & Toughness of strong material\\
    $\ell$ & $0.01$  & Damage length scale \\
    $d_1$ & $0.01$ & Relative stiffness when fully damaged  \\
    $w_1$ & $0.95$  & Damage hardening parameter \\
    \hline
    \end{tabular}\\
    \caption{Geometric and material parameters for the multi-material structures.}
    \label{tab:multiparams}
    \end{center}
    \end{table}

	\begin{figure}
		\begin{center}
			\includegraphics[width=0.9\textwidth]{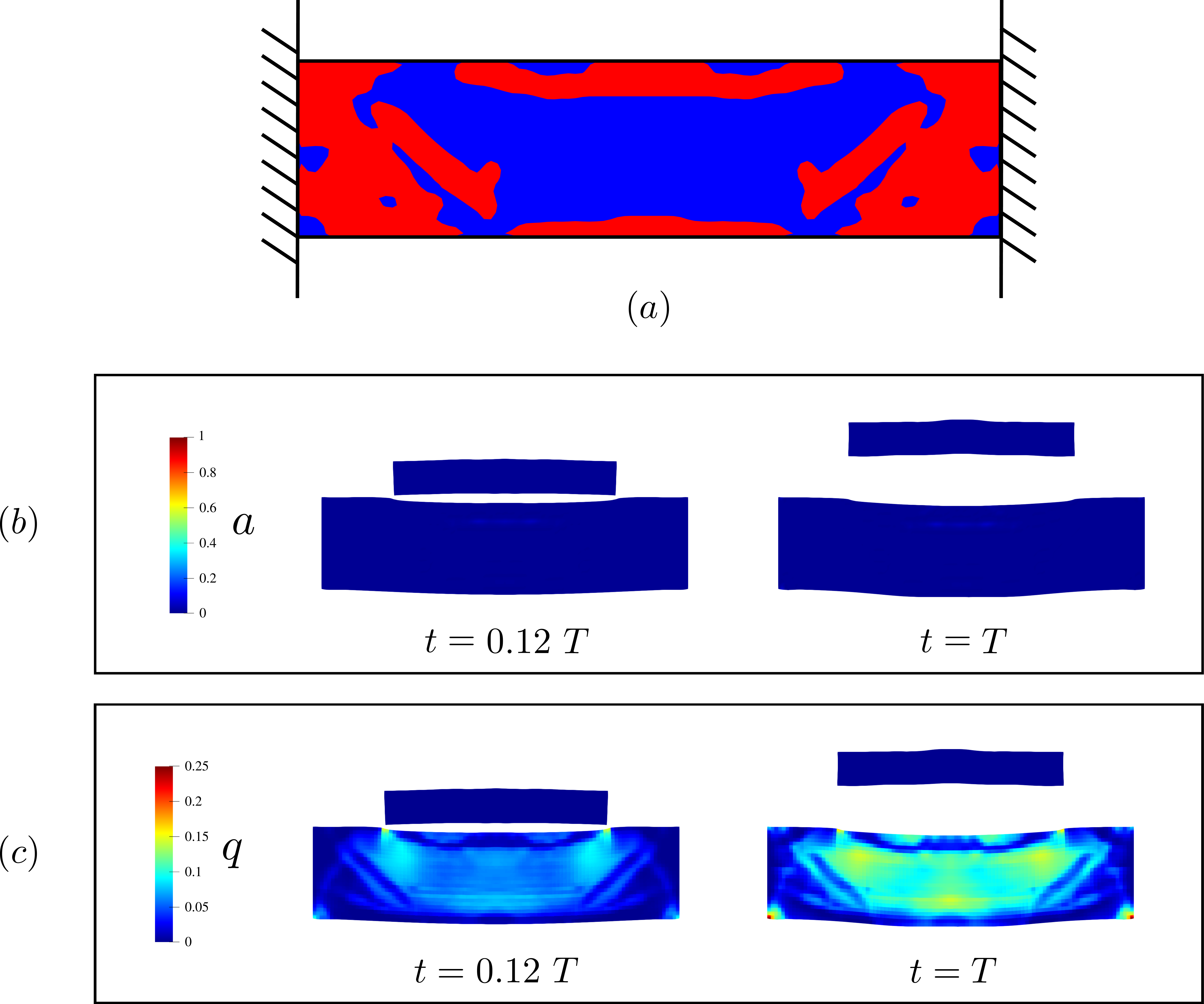}
		\end{center}
		\caption{(a) Optimal design of the multi-material structure under impact loading following contour smoothing. The red regions correspond to the strong material, and the blue regions are the tough material. (b) Damage field of this design shortly after impact (left), and long after impact (right). (c) Accumulated plasticity field of this design shortly after impact (left), and long after impact (right). }
		\label{fig:tsoptimal0}
	\end{figure}

\begin{figure}
		\begin{center}
			\includegraphics[width=1.0\textwidth]{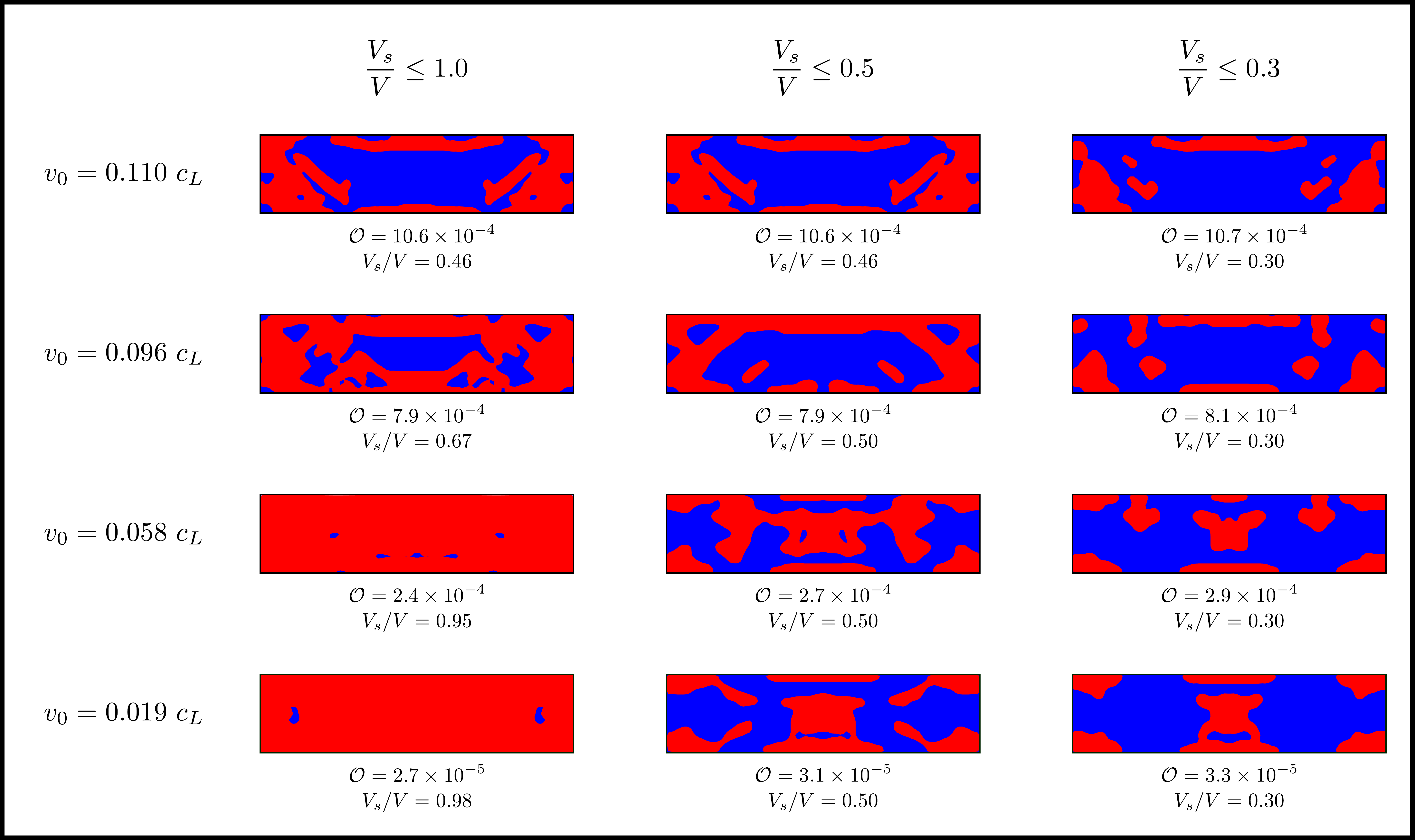}
		\end{center}
		\caption{Converged multi-material designs for impact resistance following contour smoothing. Along each row the impact velocity is constant, while along the columns the allowed amount of strong material is prescribed. $c_L$ denotes the longitudinal wave speed of the strong material. The red regions denote the strong material, while the blue regions are the tough material. Values of the objective as well as the volume fraction of strong material are shown for each design.}
		\label{fig:strongToughDesigns_vv}
	\end{figure}

    Next, we study optimal designs for varying flyer velocity and allowed volume of strong material, $V_s$. We again consider the objective in \eqref{eq:objective}. Figure~\ref{fig:strongToughDesigns_vv} shows the converged designs. For lower impact velocities, the strong material is favored. In cases where no restrictions put on the design, the converged designs are almost completely occupied by strong material. This can primarily be attributed to the stiffness difference between the strong and tough material. At $v_0 = 0.019 c_L$, there is almost no plasticity or damage, while at $v_0 = 0.058 c_L$ there is only a small amount of plasticity. However, at $v_0 = 0.110 c_L$, the converged designs have large areas of tough material, even in the case when there is no restriction placed on the amount of strong material. As discussed previously, this is to control spall which occurs at the higher impact velocities.  Additionally, strong material is used at the larger two velocities on the top surface underneath the sides of the flyer. This is to mitigate the shear-dominated plugging failure.
    
\begin{figure}
		\begin{center}
			\includegraphics[width=0.65\textwidth]{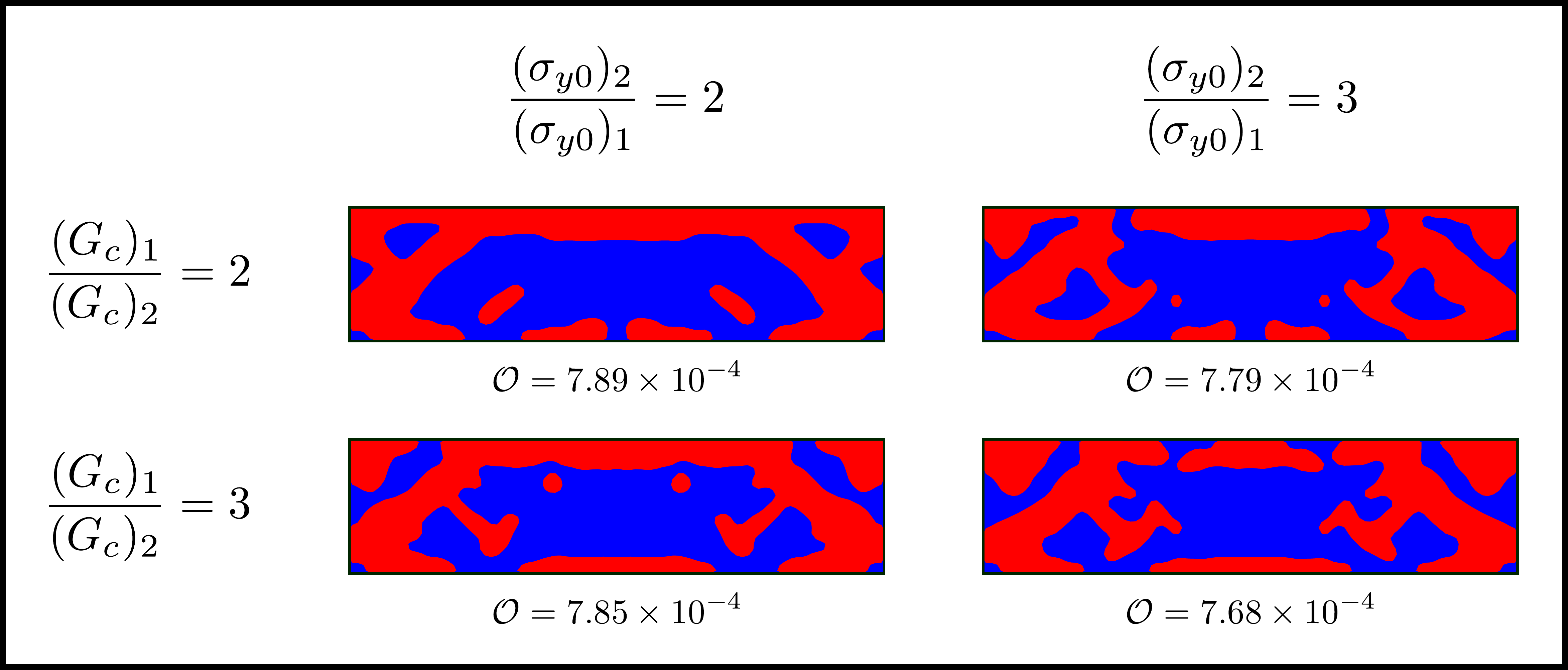}
		\end{center}
		\caption{Converged multi-material designs for impact resistance following contour smoothing. Here, we consider a constant impact velocity of $v_0 = 0.096 c_L$ and restrict $V_2/V \leq 0.5$. Along each row, the the toughness of both material are held constant, while along the columns we prescribe the yield strengths. The red regions denote the strong material, while the blue regions are the tough material. Values of the objective are shown for each design. In all of the cases, the designs saturate the constraint on allowed amount of strong material.}
		\label{fig:strongToughDesigns_ST}
	\end{figure}

	Finally, we study optimal designs for varying yield strength and toughness values. Figure~\ref{fig:strongToughDesigns_ST} shows converged designs for a constant impact velocity of $v_0 = 0.096c_L$. We constrain $V_s \leq 0.5$. The material parameters we use are identical to that of the previous study, with the following exceptions. Moving from the right to left column, the yield strength of the strong material is amplified by $50\%$ from the previous study, while moving from the top row to the bottom row has an increased toughness of the tough material by $50\%$. While the designs do vary, qualitatively they all have strong material placed near the loading site attached to struts that connect to the boundary to provide stiffness. 
	
	\section{Discussion and Conclusion} \label{sec:conclusions}

	We have developed a formulation for the optimal design of impact resistant structures. After presenting a novel method to accurately and efficiently simulating phase field damage and plasticity evolution in a transient dynamic setting, we apply gradient based optimization through the adjoint method to find optimal structures.   An important issue we address is the proper interpolation scheme for material parameters through intermediate densities. In the case of solid-void design, our formulation ensures that damage will not propagate through the void regions while preserving the natural boundary conditions at the interface. For the multi-material design, we assumed a power law interpolation for the material parameters. This implicitly penalizes intermediate densities only when either higher or lower values are clearly favorable.  These would include elastic stiffness, yield stress, and fracture toughness, where higher values are almost universally preferred. For parameters such as the damage length scale, it is unclear if a higher or lower value is favorable. However, in our study, we only consider cases where these parameters are identical for both materials. 
	
	We demonstrate these capabilities through the design of both a solid-void structure for blast loading, and a multi-material structure undergoing impact. We find that the optimal designs for the solid-void case are highly dependent on loading magnitude and duration. For the short time-scale loading, inertia plays a large role in minimizing the energy transferred to the structure. This leads to a complex trade-off between inertia and support, all while mitigating material failure. As for the multi-material structures, optimal designs use a mixture of strong and tough material when the impact velocity is high. The propogation and interaction between stress waves leads to a balance of strength and toughness throughout different parts of the domain.
	
	We now discuss possible extensions and their challenges. Not only have we simplified the contact mechanics by using asymmetric elastic elements, but we have also neglected friction and adhesion. These would be necessary to model ballistic events. It would be worthwhile to incorporate frictional contact through efficient active set methods to preserve the computational scaling~\cite{Hueber2005}. Of course, the sensitivity and adjoint formulation would need to be modified to account for this complication. Thermal effects and shock physics would be another key modeling addition. Currently, empirically derived models exist for a variety of materials which could be incorporated into this framework~\cite{Ravindran2021}.  It would also be interesting to explore other objective functionals. One might be interested in designing energy-absorbing structures that are designed to undergo plasticity and damage, rather than the objective which we chose to mitigate these. We also note that the designs we obtain depend on the location of the load. It would be straightforward to extend this work to consider multiple loading scenarios, and optimize the structure over the collective response. Finally, our simulations were done in 2D, and were readily performed on a single machine with shared memory. It would be natural to extend the implementation to a 3D settings, requiring distributed memory parallelization. 

\section*{Acknowledgement}  

The author would like to acknowledge Kaushik Bhattacharya for his advice and support throughout the project.

The financial support of the U.S. National Science Foundation through ``Collaborative Research: Optimal Design of Responsive Materials and Structures'' (DMS:2009289) and the US Army Research Laboratory thought Cooperative Agreement Number W911NF-122-0022 is gratefully acknowledged. The views and conclusions contained in this document are those of the authors and should not be interpreted as representing the official policies, either expressed or implied, of the Army Research Laboratory or the U.S. Government. The U.S. Government is authorized to reproduce and distribute reprints for Government purposes notwithstanding any copyright notation herein.

    \FloatBarrier
    
	\clearpage
	\bibliographystyle{elsarticle-num}
	\bibliography{references.bib}
	
	\newpage 
	\appendix
	
	\section{Adjoint Method for Sensitivities} \label{app:adjsens}
    We consider an objective of integral from 
    \begin{equation}
    \mathcal{O}(u, q, \varepsilon^p, a, \eta) = \int_0^T   \int_{\Omega} o(u, q, \varepsilon^p, a, \eta) \, d\Omega .
    \end{equation}
    To conduct gradient based optimization, we require the total variation of this objective with the field $\eta(x)$, which we will compute through the adjoint method. We consider adjoint fields $\xi \in \mathcal{U}$, $\gamma : \Omega \mapsto \mathbb{R}$, $\mu:\Omega \mapsto \mathbb{R}^{n \times n}$, and $b : \Omega \mapsto \mathbb{R}$ which correspond to the displacement, plastic hardening, plastic strain, and the damage field, respectively. As is standard for the adjoint method, we consider these fields as the variations in their corresponding equilibrium relations, which we add to the objective. However, for the ir-reversible damage and platicity evolution, we use the necessary Kuhn-Tucker conditions. The augmented objective is
    \begin{equation}
    \begin{aligned}
    \mathcal{O} = &\int_0^T  \int_{\Omega} \bigg\{ o + \rho \ddot{u} \cdot \xi + \pdv{W^e}{\varepsilon} \cdot  \nabla \xi - f_b \cdot \xi + \gamma \dot{q} \left[ \sigma_M - \sigma_0 - \pdv{\bar{g}^*}{\dot{q}} \right]   \\ 
    &  \qquad + \mu \cdot \left( \dot{\varepsilon}^p - \dot{q} M \right) + b \dot{a} \left[ \pdv{W^e}{a} + d^\prime \left( W^p + \int_0^t \bar{g}^* d\tau \right) - \nabla \cdot \left( \frac{G_c \ell}{2 c_w}\nabla a \right) + \frac{G_c}{4 c_w \ell }\pdv{w^a}{a} + \pdv{\bar{\psi}^*}{\dot{a}} \right]  \bigg\} \, d\Omega  dt \\
    & + \int_0^T \int_{\partial \Omega} \left( f \cdot{\xi} \right) \, dS  \, dt.
    \end{aligned}
    \end{equation}
    We then take variations with $\eta$.
    \begin{equation} \label{eq:long_adj}
    \begin{aligned}
    	\mathcal{O}_{, \eta} \delta \eta &= \int_{0}^T  \int_{\Omega} \bigg\{ \pdv{o}{\eta} + \pdv{\rho}{\eta} \ddot{u} \cdot \xi + \pdv{^2 W^e}{\varepsilon \partial \eta} \cdot \nabla \xi + b \dot{a} \left( \pdv{^2 W^e}{a \partial \eta} + \pdv{d}{a} \pdv{W^p}{\eta} + \pdv{d}{a}  \int_{0}^t \pdv{g^*}{\eta} d \tau \right) \\
    	& \quad  + \frac{1}{2 c_w}\pdv{(G_c \ell)}{\eta} \nabla(b \dot{a}) \cdot \nabla{a} + b \dot{a} \left( \frac{w^{a \prime}}{4 c_w}\pdv{(G_c/ \ell)}{\eta} + \pdv{^2 {\psi}^*}{\dot{a} \partial \eta} \right) + \gamma \dot{q} \left( \pdv{\bar{\sigma}_M}{\eta} - \pdv{\sigma_0}{\eta} - \pdv{^2 {g}^*}{\dot{q} \partial \eta}\right) \bigg\} \ \delta \eta \, d\Omega \, dt \\
    	& + \int_{0}^T \int_\Omega  \bigg\{ \pdv{o}{u} \delta_\eta u + \rho \xi \cdot \delta_\eta \ddot{u} + \left( \nabla \xi \cdot \pdv{^2 W^e}{\varepsilon \partial \varepsilon}  + b \dot{a} \pdv{^2 W^e}{a \partial \varepsilon} +  \gamma \dot{q} \pdv{\bar{\sigma}_M}{\varepsilon} - \dot{q} \mu \cdot \pdv{M}{\varepsilon} \right) \cdot \nabla \delta_\eta u  \\
    	& \quad + \left( \pdv{o}{q}  - \gamma \dot{q} \pdv{\sigma_0}{q}  + b \dot{a} d^\prime \pdv{W^p}{q}  \right) \delta_\eta q +   \left( - \gamma \dot{q} \pdv{^2\bar{g}^*}{\dot{q}^2} + \gamma  \left[ \sigma_M - \sigma_0 - \pdv{\bar{g}^*}{\dot{q}} \right] - \mu \cdot M  \right) \delta_\eta \dot{q} \\
    	& \quad + b \dot{a} d^\prime \int_0^t \left( \pdv{\bar{g}^*}{\dot{q}} \delta_\eta \dot{q}  \right) d \tau  + \mu \cdot \delta_\eta \dot{\varepsilon}^p \\
    	& \quad + \left( \pdv{o}{\varepsilon^p} + \nabla \xi \cdot  \pdv{^2 W^e}{\varepsilon \partial \varepsilon^p } + b \dot{a}  \pdv{^2 W^e}{a \partial \varepsilon}  + \gamma \dot{q}  \pdv{\bar{\sigma}_M}{\partial \varepsilon^p} - \dot{q} \mu \cdot \pdv{M}{\varepsilon^p} \right) \cdot  \delta_\eta \varepsilon^p  \\
    	& \quad + \left( \pdv{o}{a} + \pdv{^2 W^e}{a \partial \varepsilon} \cdot \nabla \xi + b \dot{a} \left[ \pdv{^2 W^e}{a^2} + \frac{G_c}{4 c_w \ell} \pdv{^2 w^a}{a^2} \right]+ b \dot{a} d^{\prime \prime} \left[  W^p + \int_0^t g^* d\tau \right) \right] \delta_\eta a  \\
    	& \quad - b \dot{a} \nabla \cdot \left( \frac{G_c \ell}{2 c_w}\nabla \delta_\eta a \right) + b D_a \delta_\eta \dot{a} + b \dot{a}  \pdv{^2\bar{\psi}^*}{\dot{a}^2} \delta_\eta \dot{a} \bigg \} \ d\Omega \ dt,
    \end{aligned}
    \end{equation}
    where
    \begin{equation}
        D_a = \pdv{W^e}{a}  + \pdv{d}{a} \left(W^p + \int_0^t g^*\, d\tau \right) - \nabla \cdot \left( \frac{G_c \ell }{2 c_w} \nabla a \right) + \frac{G_c}{4 \ell c_w }\pdv{w^a}{a}  + \pdv{\bar{\psi}^*}{\dot{a}}.
    \end{equation}
    The standard procedure would then be to integrate by parts, and enforce quiescence conditions on the adjoint variables at time $t = T$. However, for the accumulated plastic dissipation term, this is not straightforward. However, we will re-write this as
    \begin{equation}
    \begin{aligned}
        \int_0^T b \dot{a} d^\prime \int_0^t \left( \pdv{\bar{g}^*}{\dot{q}} \delta_\eta \dot{q}  \right) d \tau \ dt  &=  \int_0^T - \dv{}{t}\left [ \int_t^T b \dot{a} d^{'} \ d\tau \right]  \int_0^t \left( \pdv{\bar{g}^*}{\dot{q}} \delta_\eta \dot{q}  \right) d\tau \ dt \\
        &= - \left[ \left ( \int_t^T b \dot{a} d^{'} \ d\tau \right)  \int_0^t \left( \pdv{\bar{g}^*}{\dot{q}} \delta_\eta \dot{q}  \right) d\tau \right]_0^T + \int_0^T   \left ( \int_t^T b \dot{a} d^{'} \ d\tau \right) \pdv{\bar{g}^*}{\dot{q}} \delta_\eta \dot{q} \ dt 
    \end{aligned}
    \end{equation}
    The boundary term in the above expression is indentically zero, thus
    \begin{equation}
        \int_0^T b \dot{a} d^\prime \int_0^t \left( \pdv{\bar{g}^*}{\dot{q}} \delta_\eta \dot{q}  \right) d \tau \ dt = \int_0^T   \left ( \int_t^T b \dot{a} d^{'} \ d\tau \right) \pdv{\bar{g}^*}{\dot{q}} \delta_\eta \dot{q} \ dt 
    \end{equation}
    Using this in \ref{eq:long_adj}, we may integrate by parts. Enforcing initial quiescent conditions on the adjoint variables and localizing gives the sensitivities as
    \begin{equation} \label{eq:sensitivities_app}
    	\begin{aligned}
    		\mathcal{O}_{, \eta} \delta \eta = &\int_{0}^T  \int_{\Omega} \bigg\{ \pdv{o}{\eta} + \pdv{\rho}{\eta} \ddot{u} \cdot \xi + \pdv{^2 W^e}{\varepsilon \partial \eta} \cdot \nabla \xi + b \dot{a} \left( \pdv{^2 W^e}{a \partial \eta} + \pdv{d}{a} \pdv{W^p}{\eta} + \pdv{d}{a}  \int_{0}^t \pdv{g^*}{\eta} d \tau \right) \\
    		& + \frac{1}{2 c_w}\pdv{(G_c \ell)}{\eta} \nabla(b \dot{a}) \cdot \nabla{a} + b \dot{a} \left( \frac{w^{a \prime}}{4 c_w}\pdv{(G_c/ \ell)}{\eta} + \pdv{^2 {\psi}^*}{\dot{a} \partial \eta} \right) \\
    		& + \gamma \dot{q} \left( \pdv{\bar{\sigma}_M}{\eta} - \pdv{\sigma_0}{\eta} - \pdv{^2 {g}^*}{\dot{q} \partial \eta}\right) \bigg\} \ \delta \eta \, d\Omega \, dt,
    	\end{aligned}
    \end{equation}
    if the adjoint variables satisfy the evolution
	\begin{equation} \label{eq:adjoint_system_continuous_app}
    \begin{aligned}
    &0 = \int_{\Omega} \left[ \rho \ddot{\xi} \cdot \delta_\eta u + \pdv{o}{u} \cdot \delta_\eta u  + \left( \nabla \xi \cdot \pdv{^2 W^e}{\varepsilon \partial \varepsilon}  + b \dot{a} \pdv{^2 W^e}{a \partial \varepsilon} +  \gamma \dot{q} \pdv{\bar{\sigma}_M}{\varepsilon} - \dot{q} \mu \cdot \pdv{M}{\varepsilon} \right) \cdot \nabla \delta_\eta u \right] \ d\Omega && \forall \delta_\eta u \in \mathcal{U} \\
    & \dv{}{t} \left[ \gamma \left( \bar{\sigma}_M - \sigma_0 - \pdv{\bar{\psi}^*}{\dot{q}}\right)- \gamma \dot{q} \pdv{^2 \bar{g}^*}{\dot{q}^2} + \pdv{\bar{g}^*}{\dot{q}} \left( \int_t^T b \dot{a} d^{\prime}(a) d\tau\right) - \mu \cdot M \right]  &&\\
    & \hspace{5cm} =\pdv{o}{q} + b\dot{a} d^{\prime}(a) \pdv{W^p}{q} - \gamma \dot{q}  \pdv{\sigma_0}{q} && \text{ on } \Omega \\
    & \dv{\mu}{t} = \pdv{o}{\varepsilon^p} + \nabla \xi \cdot \pdv{^2 W^e}{\varepsilon \partial \varepsilon^p} + b \dot{a} \pdv{^2 W^e}{a \partial \varepsilon^p} + \gamma \dot{q} \pdv{\bar{\sigma}_M}{\varepsilon^p} - \dot{q} \mu \cdot \pdv{M}{\varepsilon^p} && \text{ on } \Omega \\
    &\dv{}{t} \left[ D_a b + \pdv{^2 \bar{\psi}^* }{\dot{a}^2} b \dot{a}\right] =  \pdv{o}{a} + \pdv{^2 W^e}{a \partial \varepsilon} \cdot \nabla \xi + b \dot{a} \left( \pdv{^2 W^e}{a^2} + \frac{G_c}{4 c_w \ell} \pdv{^2 w^a}{a^2}\right) &&\\
    & \hspace{5cm} + b \dot{a} d^{\prime \prime} \left(  W^p + \int_0^t g^* d\tau \right) - \nabla \cdot \left( \frac{G_c \ell}{2 c_w} \nabla(b \dot{a}) \right) && \text{ on } \Omega \\
    & \quad \xi |_{t = T} = 0, \quad  \dot{\xi} |_{t = T} = 0, \quad \gamma |_{t = T} = 0, \quad  \mu |_{t = T} = 0, \quad b |_{t = T} = 0,
    \end{aligned}
    \end{equation}

	\section{Adjoint Problem as Minimization} \label{ap:adj}
	
	It is natural to employ an augmented Lagrangian formulation to efficiently solve the adjoint problem as we have done for the forward problem. However, we will first need to write the second line of \eqref{eq:adjoint_system_continuous} as a minimization problem. Recall that this reads,
    \begin{equation}
    \begin{aligned}
     &\dv{}{t} \left[ D_a b + \pdv{^2 \bar{\psi}^* }{\dot{a}^2} b \dot{a}\right] =  \pdv{o}{a} + \pdv{^2 W^e}{a \partial \varepsilon} \cdot \nabla \xi + b \dot{a} \left( \pdv{^2 W^e}{a^2} + \frac{G_c}{4 c_w \ell} \pdv{^2 w^a}{a^2}\right) &&\\
    & \hspace{5cm} + b \dot{a} d^{\prime \prime} \left(  W^p + \int_0^t g^* d\tau \right) - \nabla \cdot \left( \frac{G_c \ell}{2 c_w} \nabla(b \dot{a}) \right) && \text{ on } \Omega \\
    \end{aligned}
    \end{equation}
    where,
    \begin{equation}
    D_a = \pdv{W^e}{a}  + \pdv{d}{a} \left(W^p + \int_0^t g^*\, d\tau \right) - \nabla \cdot \left ( \frac{G_c \ell }{2 c_w} \nabla a \right) + \frac{G_c}{4 c_w \ell} \pdv{w^a}{a}  + \pdv{\bar{\psi}^*}{\dot{a}}.
    \end{equation}
    if $\dot{a} > 0$, then $D_a = 0$. Otherwise, if $\dot{a} = 0$, then $\pdv{^2 \bar{\psi}}{\dot{a}^2} = 0$. Writing this as an implicit forward-euler discretization (as we will be solving this backwards in time) from timestep $n+1$ to $n$, gives,
    \begin{scriptsize}
    \begin{equation} \label{eq:adjoint_damage_update_raw}
    \begin{aligned}
    \frac{1}{\Delta t} \left[ \left( b^{n+1} D_a|_{t_{n+1}} + \left . \bar{\psi}^{* \prime \prime}  \right |_{t_{n + 1}} b^{n+1} \dot{a}^{n+1}  \right)  - \left .\bar{\psi}^{* \prime \prime} \right |_{t_n} b^n \dot{a}^n \right] =&   \left . \pdv{o}{a} \right|_{t_n} + \left . \pdv{^2 W^e}{a \partial \varepsilon} \right|_{t_n}  \cdot \nabla \xi^n + b^n \dot{a}^n \left . \pdv{^2 W^e}{a^2} \right|_{t_n} &&\\
    & + b^n \dot{a}^n d^{\prime \prime} \left[  W^p + \int_0^t g^* d\tau \right]_{t_n} - \nabla \cdot \left( \frac{G_c \ell}{2 c_w} \nabla(b^n \dot{a}^n) \right) \quad &&  \text{ on }\Omega_{\dot{a}_n > 0} \\
    \frac{1}{\Delta t} \left[ b^n D_a |_{t_n} - \left( b^{n+1} D_a|_{t_{n+1}} + \left . \bar{\psi}^{* \prime \prime}  \right |_{t_{n + 1}} b^{n+1} \dot{a}^{n+1}  \right) \right] =&   \left . \pdv{o}{a} \right|_{t_n} + \left . \pdv{^2 W^e}{a \partial \varepsilon} \right|_{t_n}  \cdot \nabla \xi^n \qquad && \text{ on } \Omega_{\dot{a}_n = 0} 
    \end{aligned}
    \end{equation}
    \end{scriptsize}
    If we define,
    \begin{equation}
    z^{n} = \dot{a}^n b^n ,
    \end{equation}
    we may write the first line of \eqref{eq:adjoint_damage_update_raw} as a minimization problem
    \begin{equation}
    \begin{aligned}
    \inf_{z = 0 \: \text{on} \: \Omega_{\dot{a}_n = 0}}  I[z] = \int_{\Omega} \Bigg \{ &\frac{1}{2 \Delta t \bar{\psi}^{* \prime \prime} |_{t_n}} \left[ \left( b^{n+1} D_a|_{t_{n+1}} + \left .\bar{\psi}^{* \prime \prime}  \right |_{t_{n + 1}} b^{n+1} \dot{a}^{n+1}  \right) - \left . \bar{\psi}^{* \prime \prime}  \right |_{t_n} z  \right]^2 \\
    & + \left(  \left . \pdv{o}{a} \right|_{t_n} + \left . \pdv{^2 W^e}{a \partial \varepsilon} \right|_{t_n}  \cdot \nabla \xi^n \right) z  \\
    & + \left( \left . \pdv{^2 W^e}{a^2} \right|_{t_n} +  d^{\prime \prime} \left[  W^p + \int_0^t g^* d\tau \right]_{t_n} \right) \frac{z^2}{2} + \frac{G_c \ell}{4 c_w} \abs{\nabla z}^2 \Bigg \} \ d\Omega.
    \end{aligned}
    \end{equation}
    We will now introduce another augmented Lagrangian with an auxiliary field $\zeta \in L^2(\Omega)$, and enforce $\zeta = z$  through the Lagrange multiplier field $\chi \in L^2(\Omega)$ and penalty factor $r$. Thus, the previous minimization is equivalent to finding the saddle point of,
    \begin{equation}
    \begin{aligned}
    \hat{\mathcal{L}}(z, \zeta, \chi) =  \int_{\Omega} \Bigg \{ &\frac{1}{2 \Delta t \bar{\psi}^{* \prime \prime} |_{t_n}} \left[ \left( b^{n+1} D_a|_{t_{n+1}} + \left .\bar{\psi}^{* \prime \prime}  \right |_{t_{n + 1}} b^{n+1} \dot{a}^{n+1}  \right) - \left . \bar{\psi}^{* \prime \prime}  \right |_{t_n} z  \right]^2 \\
    & + \left(  \left . \pdv{o}{a} \right|_{t_n} + \left . \pdv{^2 W^e}{a \partial \varepsilon} \right|_{t_n}  \cdot \nabla \xi^n \right) z  + \left( \left . \pdv{^2 W^e}{a^2} \right|_{t_n} +  d^{\prime \prime} \left[  W^p + \int_0^t g^* d\tau \right]_{t_n} \right) \frac{z^2}{2}\\
    &  + \frac{G_c \ell}{4 c_w} \abs{\nabla z}^2 + \chi \left( z - \zeta \right) + \frac{r}{2} (z - \zeta)^2 \Bigg \} \ d\Omega
    \end{aligned}
    \end{equation}
    subject to the constraints that $\zeta = 0$ on $\Omega_{\dot{a} = 0}$. Then, conditions for stationarity are
    \begin{equation}
    \begin{aligned}
    &0 = \int_{\Omega} \left[ \left( r(z - \zeta^n) + \chi \right) \delta z + \frac{G_c \ell}{2 c_w} \nabla z \cdot \nabla \delta z \right] \, d \Omega && \hspace{-1cm}\forall \delta z \in \mathcal{A} \\
    &0 = \int_{\Omega} \left(  z - \zeta^n \right) \delta \chi \, d \Omega && \hspace{-1cm} \forall \delta \chi \in L^2(\Omega) \\ 
    \end{aligned}
    \end{equation}
    and
    \begin{equation}
    \begin{cases}
    \frac{1}{\Delta t} \left[ \left( b^{n+1} D_a|_{t_{n+1}} + \left . \bar{\psi}^{* \prime \prime}  \right |_{t_{n + 1}} \zeta^{n+1}  \right)  - \left .\bar{\psi}^{* \prime \prime} \right |_{t_n} \zeta^n \right] =   \left . \pdv{o}{a} \right|_{t_n} + \left . \pdv{^2 W^e}{a \partial \varepsilon} \right|_{t_n}  \cdot \nabla \xi^n + \zeta^n \left . \pdv{^2 W^e}{a^2} \right|_{t_n} &\\
    \hspace{8cm} + \zeta^n d^{\prime \prime} \left[  W^p + \int_0^t g^* d\tau \right]_{t_n} - r(z - \zeta^n) - \chi \qquad   &\text{ on }\Omega_{\dot{a}_n > 0} \\
    \zeta^n = 0 \qquad  &\text{ on } \Omega_{\dot{a}_n = 0} \\
    \frac{1}{\Delta t} \left[ \left( b^{n+1} D_a|_{t_{n+1}} + \left . \bar{\psi}^{* \prime \prime}  \right |_{t_{n + 1}} \zeta^{n+1}  \right) - D_a |_{t_n} b^n \right] =   \left . \pdv{o}{a} \right|_{t_n} + \left . \pdv{^2 W^e}{a \partial \varepsilon} \right|_{t_n}  \cdot \nabla \xi^n - r z - \chi  \qquad  & \text{ on } \Omega_{\dot{a}_n = 0} 
    \end{cases}
    \end{equation}
    Using $\zeta^n = \dot{a}^n b^n$ gives
    \begin{equation}
    \begin{aligned}
    &0 = \int_{\Omega} \left[ \left( r(z - \dot{a}^n b^n) + \chi \right) \delta z + \frac{G_c \ell}{2 c_w} \nabla z \cdot \nabla \delta z \right] \, d \Omega && \forall \delta z \in \mathcal{A} \\
    &0 = \int_{\Omega} \left(  z -\dot{a}^n b^n \right) \delta \chi \, d \Omega &&  \forall \delta \chi \in L^2(\Omega) \\ 
    &\frac{1}{\Delta t} \left[ \left( b^{n+1} D_a|_{t_{n+1}} + \left . \bar{\psi}^{* \prime \prime}  \right |_{t_{n + 1}} \dot{a}^{n+1} b^{n+1}  \right)  - \left .\bar{\psi}^{* \prime \prime} \right |_{t_n} \dot{a}^n b^n - b^{n} D_a|_{t_{n}}\right] =  \\
    & \qquad \left . \pdv{o}{a} \right|_{t_n} + \left . \pdv{^2 W^e}{a \partial \varepsilon} \right|_{t_n}  \cdot \nabla \xi^n + \dot{a}^n b^n \left . \pdv{^2 W^e}{a^2} \right|_{t_n} + \dot{a}^n b^n d^{\prime \prime}   \left[  W^p + \int_0^t g^* d\tau \right]_{t_n} - r(z - \dot{a}^n b^n) - \chi \qquad  &&\text{ on } \Omega 
    \end{aligned}
    \end{equation}
	
\end{document}